\documentclass{article}

\usepackage{arxiv}

\usepackage[utf8]{inputenc} 
\usepackage[T1]{fontenc}    
\usepackage{url}            
\usepackage{booktabs}       
\usepackage{amsfonts}       
\usepackage{nicefrac}       
\usepackage{microtype}      
\usepackage{lipsum}
\usepackage{amsmath}
\usepackage{amsthm}
\usepackage{amssymb}
\usepackage{bm}
\usepackage{diagbox}
\usepackage{caption}
\usepackage{subcaption}
\usepackage{graphicx}
\usepackage[ruled,vlined]{algorithm2e}
\include{pythonlisting}
\usepackage{xcolor}
\usepackage{mathtools}
\usepackage{enumitem}
\usepackage{cite}
\usepackage{hyperref}
\hypersetup{
    colorlinks=true,
    citecolor=black,
    linkcolor=black,
    filecolor=black,      
    urlcolor=black,
}

\newcommand{\R}{\mathbb{R}}

\theoremstyle{definition}

\numberwithin{equation}{section}

\title{Deep learning of free boundary and Stefan problems}

\author{
  Sifan Wang \\
  Graduate Group in Applied Mathematics \\
  and Computational Science \\
  University of Pennsylvania\\
  Philadelphia, PA 19104 \\
  \texttt{sifanw@sas.upenn.edu} \\
   \And
  Paris Perdikaris \\
  Department of Mechanichal Engineering \\
  and Applied Mechanics\\
  University of Pennsylvania\\
  Philadelphia, PA 19104 \\
  \texttt{pgp@seas.upenn.edu} \\
}

\begin{document}
\maketitle

\begin{abstract}
Free boundary problems appear naturally in numerous areas of mathematics, science and engineering. These
problems present a great computational challenge because they necessitate numerical methods that can yield an accurate approximation of free boundaries and complex dynamic interfaces. In this work, we propose a multi-network model based on physics-informed neural networks to tackle a general class of forward and inverse free boundary problems called Stefan problems. Specifically, we approximate the unknown solution as well as any moving boundaries by two deep neural networks. Besides, we formulate a new type of inverse Stefan problems that aim to reconstruct the solution and free boundaries directly from sparse and noisy measurements. We demonstrate the effectiveness of our approach in a series of benchmarks spanning different types of Stefan problems, and illustrate how the proposed framework can accurately recover solutions of  partial differential equations with moving boundaries and dynamic interfaces.
All code and data accompanying this manuscript are publicly available at \url{https://github.com/PredictiveIntelligenceLab/DeepStefan}.
\end{abstract}


\section{Introduction}
\label{sec:Intro}


Free boundary problems define a class of mathematical problems involving the solution of partial differential equations in a domain, a part of whose boundary is a-priori unknown and has to be determined as part of the problem's solution \cite{chen2015free}. Such problems are ubiquitous in science and engineering as they can naturally model continuum systems with phase transitions, moving obstacles, multi-phase dynamics, competition for resources, etc., with classical applications ranging from tumor growth modeling \cite{rutter2017mathematical, friedman1999stefan, liu2018analysis, chaplain1991mathematical}, to chemical vapor deposition \cite{friedman2000free},  crystal growth and solidification in materials \cite{madejski1976solidification}, semi-conductor design \cite{friedman1995stefan}, and option pricing \cite{pham1997optimal}.
The mathematical origin of free boundary problems traces back to the famous {\em Stefan problem} that describes the joint evolution of a liquid and a solid phase related to heat transfer \cite{stefan1891theorie}; a question also considered in 1831 by Lam\'e and Clapeyron in relation to the problems of ice formation in the polar seas \cite{chen2015free, tarzia2000bibliography}. Since their initial conception nearly two centuries ago, free boundary problems now define a well-studied area in applied mathematics both in terms of theory \cite{duvant2012inequalities}, numerical methods \cite{crowley1978numerical, ronquist1987legendre, chen1997simple}, and applications across a wide range of problems in science and engineering \cite{friedman2000free}.

Free boundary and Stefan problems can be generally classified into two distinct categories, namely direct and inverse problems. In its classical form, the direct Stefan problem involves finding the temperature distribution in a domain undergoing a phase transition, as well as determining the position of the free boundary defining the transition interface. The latter is identifiable  thanks to a number of free boundary conditions -- the so called Stefan conditions -- that are derived from energy conservation or other physical laws 
governing the phase transition interface \cite{friedman2000free}. Variations of the Stefan problem also include cases in which one or several variables must be determined in different domains of the space, or space–time, whereas in each domain a different set of state laws and evolution equations may be specified \cite{chen2015free, friedman2000free}. On the other hand, the definition inverse free boundary problems is still somewhat open ended \cite{gol2012inverse}. For example, in its classical form, the inverse Stefan problem is focused on either inferring unknown thermophysical parameters, or determining an initial condition or a boundary condition with the help of additional information such as the position of the moving boundary interface or measurements of temperature distribution at some time instants. These inverse Stefan-type problems are always expected to be ill-posed and thus some regularization methods are needed to obtain stable numerical solutions, where ``ill-posed" means that the solution does not continuously depend on the input  data.

Different numerical methods have been developed to solve various direct and inverse Stefan problems. Depending on how they treat the free boundaries and phase transition interfaces, these methods can be generally divided into explicit and implicit methods, respectively, see \cite{furzeland1980comparative,javierre2006comparison}. One representative example of explicit methods is the front tracking method \cite{womble1989front, juric1996front, tryggvason2001front, unverdi1992front}, which aims to explicitly track the moving interfaces in a time-varying domain . If accurately resolved, these methods can yield highly accurate moving boundaries, however they have difficulty in handling moving boundaries that may develop irregular shapes (e.g. sharp cusps, double backs, or even disappear), especially in high dimensions. In contrast to explicit methods, the computational grid in implicit methods is usually fixed and the position of the interface is obtained indirectly from some field defined over the whole domain. The key advantage here is that the moving interfaces need not be tracked explicitly, hence avoiding the key limitations of explicit approaches. Examples of implicit methods include enthalpy methods \cite{voller1987implicit, date1992novel}, phase-field \cite{fix1982phase, mackenzie2002moving} and level-set methods \cite{chen1997simple, zheng2000adaptive, osher2004level, olsson2005conservative}). The enthalpy formulation is one of the most popular fixed domain methods for solving Stefan problems by reformulating the governing equations in terms of enthalpy such that any free boundary conditions are automatically satisfied. Despite of its obvious advantages and favorable computational cost, this approach is only applicable to a relatively narrow range of problems for which the enthalpy reformulation can be derived in closed form. Phase-field methods retain a fixed grid approach and define a phase function $\phi$, which is constant within each phase and varies smoothly between these values in the interface regions. However, the interface position cannot be defined exactly because an artificial interface thickness must be added to the model to allow smooth variation. Finally, level set methods aim to capture the position of the interface as represented by the zero-level set of a continuous function. Some comparative analyses \cite{javierre2006comparison} lead to the conclusion that methods based on the level set are quite promising and general. However, such methods are typically not volume-preserving and therefore prone to numerical instabilities. More recently, methods based on deep learning have been proposed to  solve high-dimensional free boundary problems \cite{sirignano2018dgm}, however they rely on the assumption that the free boundary is known and fixed.

All aforementioned methods have their own advantages and limitations and have been proven effective for certain classes of direct Stefan problems. Some of these methods can also serve as building blocks for solving inverse Stefan problems, which typically require the repeated evaluation of direct solvers within an adjoint optimization loop \cite{jochum1980numerical,colton1984numerical,slota2007direct,johansson2011method}. However, in most cases these methods have to be specifically tailored to a given problem setup and cannot be easily adapted to build a general framework for seamlessly tackling direct and inverse free boundary problems involving various types of partial differential equations. Moreover, to the best of our knowledge, none of these existing techniques can be readily applied to data assimilation problems in which no initial or boundary conditions may be known, but instead a finite number of sparse measurements inside the domain of interest is available. 

In this work, we expand on the recently developed framework of physics-informed neural networks \cite{raissi2019physics} to tackle a wide range of direct and inverse free boundary and Stefan problems.
Our specific contributions can be summarized in the following points:
\begin{itemize}[leftmargin=*]
    \item We propose a unified computational framework that can be applied to effectively tackle all kinds of Stefan and free boundary problems. 
    \item We demonstrate the effectiveness of our method in solving variable types of direct and inverse Stefan problems. 
    \item  We also propose and solve a new type of data-driven inverse Stefan problem that does not require initial and boundary conditions, but instead utilizes sparse measurements of the dependent variables in a given domain.
\end{itemize}
Taken all together, the proposed methods present a new general framework for tackling free boundary and Stefan problems that commonly appear in biology, physics, engineering, and finance.

This paper is organized as follows. In section \ref{sec: overview}, we present a brief overview of the physics-informed neural networks (PINNs) \cite{raissi2019physics}. In section \ref{sec: problem_formulation}, we outline the proposed deep learning framework and demonstrate its effectiveness through the lens of three representative case studies, including one-dimensional one-phase Stefan problems, one-dimensional two-phase Stefan problems, and two-dimensional one-phase Stefan problems. Finally, we summarize our findings and discuss promising future directions in section \ref{sec: summary}.

\section{A brief overview on physics-informed neural networks}
\label{sec: overview}

Under the framework of physics-informed neural networks (PINNs) \cite{raissi2019physics}, we consider partial differential equations (PDEs) of the following general form

\begin{align}
\label{eq:PDE}
    \begin{split}
     &\bm{u}_t + \mathcal{N}_{\bm{x}}[\bm{u}] = 0, \ \  \bm{x} \in \Omega, t \in [0, T] \\
     &\bm{u}(\bm{x}, 0) = h(\bm{x}), \ \ \bm{x} \in \Omega \\
     &\bm{u}(\bm{x}, t) = g(\bm{x}, t), \ \ t \in [0, T], \ \  \bm{x} \in \partial \Omega,
    \end{split}
\end{align} 
where $\bm{x}\in\mathbb{R}^{d}$ and $t\in[0,T]$ are spatial and temporal coordinates, $\Omega$ denotes a bounded domain in $\mathbb{R}^{d}$ with boundaries $\partial \Omega$, $T>0$, and  $\mathcal{N}_{\bm{x}}$ is a nonlinear differential operator. In addition, $u(\bm{x},t) : \Omega \times [0,T] \rightarrow \R$ describes the unknown physical law governed by the above equations \ref{eq:PDE}. 

Following the original work of Raissi {\em et. al.} \cite{raissi2019physics} we assume that the latent solution $u(\bm{x},t)$ can be approximated by a deep neural network $u_{\bm{\theta}}(\bm{x},t)$ with parameters $\bm{\theta}$. Then the boundary residual $\mathcal{L}_{u_b}(\bm{\theta})$, the initial residual $\mathcal{L}_{u_0}(\bm{\theta})$, and the PDE residual $\mathcal{L}_r(\bm{\theta})$ can be defined as follows
\begin{align}
\label{eq: loss_r}
    &\mathcal{L}_r(\bm{\theta}) = \frac{1}{N_r} \sum_{i=1}^{N_r} [\bm{r}(\bm{x}_r^i, t_r^i)]^2 \\
    \label{eq: loss_ub}
     &\mathcal{L}_{u_b}(\bm{\theta}) =\frac{1}{N_b} \sum_{i=1}^{N_b}[\bm{u}_{\bm{\theta}}(\bm{x}_b^i, t_b^i) - g(\bm{x}_b^i, t_b^i)]^2, \\
     \label{eq: loss_u0}
    &\mathcal{L}_{u_0}(\bm{\theta}) = \frac{1}{N_0} \sum_{i=1}^{N_0}[\bm{u}_{\bm{\theta}}(\bm{x}_0^i, 0) - h(\bm{x}_0^i))]^2,
\end{align}
where,
\begin{align}
    r_{\theta}(\bm{x}, t):=\frac{\partial}{\partial t} u_{\theta}(\bm{x}, t)+\mathcal{N}_{x}\left[u_{\theta}(\bm{x}, t)\right].
\end{align}
The parameters of the neural network $u_{\bm{\theta}}(\bm{x},t)$ can be estimated by minimizing the mean square loss function of the form
\begin{align}
    \label{eq: loss}
    \mathcal{L}(\bm{\theta}) = \lambda_r \mathcal{L}_r(\bm{\theta}) + \lambda_b \mathcal{L}_{u_b}(\bm{\theta}) + \lambda_{u_0} \mathcal{L}_{u_0}(\bm{\theta}).
\end{align}
Note that the parameters $\left\{ \lambda_r,  \lambda_b,  \lambda_0 \right\}$ denote  weight coefficients in the loss function, and can effectively assign a different learning rate to each individual loss term. These weights may be user-specified or tuned manually or automatically by utilizing the back-propagated gradient statistics during network training \cite{wang2020understanding}.

Thanks to the approximation capabilities of neural networks,  physics-informed neural networks have already led to a series of noticeable results across a range of problems in computational science and engineering, including stochastic differential equations \cite{han2018solving, zhang2020learning}, fractional differential equations \cite{pang2019fpinns}, uncertainty quantification \cite{yang2019adversarial,zhu2019physics, yang2018physics}, bio-engineering    \cite{sahli2020physics, kissas2020machine}, meta-material design \cite{fang2019deep,liu2019multi,chen2020physics}, fluids mechanics \cite{raissi2020hidden,sun2020surrogate,raissi2019deep1,raissi2019deep2,jin2020nsfnets}, and beyond \cite{tartakovsky2018learning,lu2019deeponet,tartakovsky2020physics}. All cases considered so far in the literature pertain to problems where the domain boundaries are precisely known and no interfacial phenomena take place. Within these well-defined domains, it is natural to employ neural networks to approximate the latent solution of the PDE. However, in free boundary and Stefan problems, not only is the latent solution unknown, but parts of the domain's boundaries may also be unknown. So it is also natural that such free boundaries $s(t)$ can be parametrized by another neural network, as discussed in the following sections.

\section{Case studies}
\label{sec: problem_formulation}

The goal of the following subsections is to outline the key technical details that define the proposed deep learning framework, and illustrate its effectiveness across a range of free boundary and Stefan problems, including one or two dimensional Stefan problems as well as one-phase or two-phase Stefan problems. In each class of problems, we will consider both direct and inverse problems. In addition, we formulate a new type of inverse Stefan problems which has great potential for meeting various demands in the data-driven setting. 

Throughout all case studies we will use a universal algorithmic setup. Fully-connected neural networks will be used to approximate the latent functions representing PDE solutions and unknown boundaries. We employ hyperbolic tangent activation functions, stochastic gradient updates using a mini-batch size of 128 data-points, and train the networks using the Adam optimizer with default settings \cite{kingma2014adam}. Moreover, all networks are initialized using the Glorot scheme \cite{glorot2010understanding}, and no additional regularization techniques are employed (e.g., dropout, $L^1$/$L^2$ penalties, etc.). These hyper-parameters settings are summarized in table \ref{tab: hyper-parameters}, and are kept fixed throughout all numerical studies. All code and data accompanying this manuscript are publicly available at \url{https://github.com/PredictiveIntelligenceLab/DeepStefan}. The training time for performing $40,000$ stochastic gradient iterations is roughly 5 minutes for all numerical experiments on a Lenovo X1 Carbon ThinkPad laptop with an Intel Core i5 2.3GHz processor and 8Gb of RAM memory.

\begin{table}[]
   \renewcommand*{\arraystretch}{1.4}
    \centering
   \begin{tabular}{|c|c|}
\hline
Hyperparameter          & Value  \\ \hline
learning rate           & $10^{-3}$   \\ \hline
iterations              & 40,000 \\ \hline
batch-size              & 128    \\ \hline
number of hidden layers & 3      \\ \hline
number of hidden units  & 100    \\ \hline
activation function     & Tanh   \\ \hline
Initialization          & Xavier \\ \hline
Optimizer               & Adam \\ \hline
\end{tabular}
    \caption{Hyper-parameter settings employed throughout all numerical experiments presented in this work.}
    \label{tab: hyper-parameters}
\end{table}

\subsection{One-dimensional one-phase Stefan problems}
\label{sec: Stefan_1D1P}

To begin our presentation, and to give a clear mathematical formulation of the proposed methodology, we first restrict ourselves to defining the classical one-dimensional one-phase Stefan problem \cite{furzeland1980comparative}, which considers a semi-infinite solid, for instance a thin block of ice occupying $0 \leq x < \infty$ at the solidification temperature $u=0$. At the fixed boundary of the thin block of ice $x=0$ there
could be different types of ``flux functions" $g(t)$.
For any $t > 0$, the region $0 \leq x < \infty$ will consist of solid and liquid phases with the liquid
phase occupying the region $0 \leq x < s(t)$ and the solid phase the region $s(t) < x < \infty$. This ice melting problem can be formulated mathematically as follows.

Let  $\Omega = (0, s(t)) \times (0, T)$, where $s(t)$ denotes an unknown  moving boundary. In the domain $\Omega$, the unknown temperature distribution $u(x,t)$ satisfies the one-dimensional heat equation
\begin{align}
      \label{eq: Stefan_1D1P_eq}
     \frac{\partial u}{\partial t}  - \frac{\partial^2 u}{\partial x^2}  =0,  \quad (x,t) \in  (0,s(t)) \times (0, T),
\end{align}
subject to initial and Neumann boundary conditions 
\begin{align}
      \label{eq: Stefan_1D1P_ic}
    & u(x, 0) =  u_0(x),  \quad x \in [0, s(0)] \\
     \label{eq: Stefan_1D1P_uNc}
    &  \frac{\partial u}{\partial x}(0, t) = g(t)  , \quad t \in [0, T].
\end{align}
The Dirichlet and Neumann boundary conditions on the moving boundary $x= s(t)$ are given by
\begin{align}
 \label{eq: Stefan_1D1P_Sic}
    &s(0) = s_0 \\
  \label{eq: Stefan_1D1P_Sbc}
    & u(s(t), t) = h_1(t), \quad t \in [0, T] \\
    \label{eq: Stefan_1D1P_SNc}
    & \frac{\partial u}{\partial x}(s(t), t) =  h_2(t), \quad t \in [0, T].
\end{align}
Equation \ref{eq: Stefan_1D1P_Sic} describes the initial position of the melting interface, while equation \ref{eq: Stefan_1D1P_Sbc} describes the  
equilibrium phase-change temperature as both solid and liquid phases can stay together in thermodynamic equilibrium. For instance, $h_1(t) = 0$ is the freezing temperature.


For direct Stefan problems, $u_0, g(t), h_1(t), h_2(t), s_0$ are assumed to be known, and our goal is to infer the unknown temperature solution $u(x,t) : \Omega \rightarrow \R$, as well as determine the moving boundary $s(t)$ satisfying equation \ref{eq: Stefan_1D1P_eq} - \ref{eq: Stefan_1D1P_Sic}.

In contrast to direct Stefan problem, there are many types of inverse Stefan problems that can be formulated. One classical inverse Stefan problem -- coined here as Inverse Type I -- is the  moving boundary design problem \cite{reemtsen1984method,slota2008solving,johansson2011method}. In this type of problem, we wish to find the temperature distribution $u(x, t)$ satisfying \ref{eq: Stefan_1D1P_eq} -\ref{eq: Stefan_1D1P_ic}, \ref{eq: Stefan_1D1P_Sbc} - \ref{eq: Stefan_1D1P_SNc}, as well as reconstruct the Dirichlet and Neumann boundary conditions when the position of moving boundary $s(t)$ is known in advance.

However, in many realistic cases, it is impractical to track the  position of the phase-change boundary. Conversely, it is relatively easier to obtain some measurements of temperature in the given domain. This motivates us to propose a new class of inverse Stefan problems. That is, instead of considering any initial or boundary conditions, some temperature measurements in the domain $\Omega$ are given. We then seek to find the temperature solution $u(x,t)$ and identify the unknown position of the moving boundary, satisfying \ref{eq: Stefan_1D1P_eq} and \ref{eq: Stefan_1D1P_Sbc} - \ref{eq: Stefan_1D1P_Sic}. We will refer to this class of Stefan problems as Inverse Type II. Note that this type of inverse Stefan problem is different from the formulation outlined in Benard {\em et. al.} \cite{benard1992inverse}, in which discrete measurements of temperatures and fluxes can be collected at the fixed part of the boundary. In our case, neither do we need measurements of heat fluxes nor have restrictions on the location of the experimental data. 
 
Table \ref{tab: Stefan_1D1P} summarizes the different classes of one-dimensional one-phase Stefan problems that we will be presenting as illustrative cases for the application of the proposed algorithms.
 
 \begin{table}
    \centering
     \renewcommand*{\arraystretch}{1.8}
   \begin{tabular}{|c|c|c|c|}
\hline
                & Equations  & Observed & Latent \\ \hline
Direct          &  \ref{eq: Stefan_1D1P_eq} - \ref{eq: Stefan_1D1P_Sic}         &  $u_0(x), g(t), h_1(t),h_2(t), s_0$ & $u(x,t), s(t)$       \\ \hline
Inverse Type I  &  \ref{eq: Stefan_1D1P_eq} \ref{eq: Stefan_1D1P_ic} \ref{eq: Stefan_1D1P_Sbc} \ref{eq: Stefan_1D1P_SNc}             & $s(t), u_0(x),h_1(t),h_2(t)$ & $u(x,t), u_x(x,0)$     \\ \hline
Inverse Type II & \ref{eq: Stefan_1D1P_eq}  \ref{eq: Stefan_1D1P_Sbc} - \ref{eq: Stefan_1D1P_Sic}        & $\{(x^j_{\text{data}}, t^j_{\text{data}}, u^j)\}_{j=1}^M, h_1(t),h_2(t),s_0$    &  $u(x,t), s(t)$    \\ \hline
\end{tabular}
    \caption{{\em One-dimensional one-phase Stefan problem:} Summary of conditions and the objective of each type of Stefan problem formulated in section \ref{sec: Stefan_1D1P}.}
    \label{tab: Stefan_1D1P}
\end{table}


\subsubsection{Direct one-dimensional one-phase Stefan problems}
\label{sec: Stefan_1D1P_direct}

To illustrate the proposed workflow, let us focus on a concrete benchmark example for which an analytical solution can be derived \cite{johansson2011method}. To this end, consider the direct one-dimensional one-phase Stefan problem described in section \ref{sec: problem_formulation} with an artificial domain $\Omega^* = [0,1] \times [0,1] \supset \Omega=(0, s(t)) \times (0,1)$. The known variables in \ref{eq: Stefan_1D1P_eq} - \ref{eq: Stefan_1D1P_Sic} are given by
\begin{align}
     &s_0 = 2 - \sqrt{3}, \quad u_0 = -\frac{x^2}{2} + 2x - \frac{1}{2}, \quad x \in [0, s(0)], \\
     &g(t) = 2, \quad h_1(t) = 0, \quad h_2(t) = \sqrt{3 - 2t}, \quad t\in [0,1],
\end{align}
Then, the exact solution and the moving boundary are given by
\begin{align}
    \label{eq: Stefan_1D1P_sol}
    &u(x,t) = -\frac{x^2}{2} + 2x - \frac{1}{2} - t, \quad (x,t) \in  [0,s(t)] \times [0, 1] \\
    \label{eq: Stefan_1D1P_mb_sol}
    & s(t) = 2 - \sqrt{3 -2t}, \quad t \in [0, 1].
\end{align}

Recall that our goal is to find the temperature solution $u(x,t)$ and the moving boundary $s(t)$ satisfying equations \ref{eq: Stefan_1D1P_eq} and \ref{eq: Stefan_1D1P_Sic}. To this end, we use two independent deep neural networks $u_{\bm \theta}(x,t)$ and $s_{\bm \beta}(t)$ to approximate the latent solution $u(x,t)$ and $s(t)$, where ${\bm \theta}$ and ${\bm \beta}$ are two separate parameter spaces, respectively. We expect that a physics-informed neural network model $u_{\bm \theta}(x,t)$ with parametrized boundary $s_{\bm \beta}(t)$ can now be trained to approximate the latent solution $u(x,t)$ as well as the unknown free boundary $s(t)$ by minimizing the following composite loss
\begin{align}
    \label{eq: loss_Stefan_1D1P_direct}
    \mathcal{L}({\bm \theta},\bm{\beta} ) =   \mathcal{L}_r({\bm \theta}) +  \mathcal{L}_{u_0}({\bm \theta}) +  \mathcal{L}_{u_{Nc}}({\bm \theta})  +  \mathcal{L}_{s_{bc}}({\bm \theta},\bm{\beta}) +  \mathcal{L}_{s_{Nc}}({\bm \theta},\bm{\beta}) + \mathcal{L}_{s_0}(\bm{\beta}),
\end{align}
where,
\begin{align}
    \label{eq: loss_Stefan_1D1P_L_r}
    & \mathcal{L}_r({\bm \theta} ) = \frac{1}{N} \sum_{i=1}^{N}|\frac{\partial u_{\bm \theta}}{\partial t}(x^i, t^i) - \frac{\partial^2 u_\theta}{\partial x^2}(x^i, t^i)|^2  \\
     \label{eq: loss_Stefan_1D1P_L_u0}
    &\mathcal{L}_{u_0}({\bm \theta}) = \frac{1}{N} \sum_{i=1}^{N} |u_{\bm \theta}(x^i,0) - u_0(x^i)|^2 \\
     \label{eq: loss_Stefan_1D1P_L_uNc}
    & \mathcal{L}_{u_{Nc}}({\bm \theta}) = \frac{1}{N} \sum_{i=1}^N|\frac{\partial u_{\bm \theta}}{\partial x}(0, t^i) - g(t^i)|^2  \\
     \label{eq: loss_Stefan_1D1P_L_sbc}
    & \mathcal{L}_{s_{bc}}({\bm \theta},\bm{\beta}) = \frac{1}{N} \sum_{i=1}^{N} |u_{\bm \theta}(s_{\bm \beta}(t^i), t^i)|^2 \\
     \label{eq: loss_Stefan_1D1P_L_sNc}
    & \mathcal{L}_{s_{Nc}}({\bm \theta},\bm{\beta}) = \frac{1}{N} \sum_{i=1}^{N} |\frac{\partial u_{\bm \theta}}{\partial x}(s_{\bm \beta}(t^i), t^i) - h(t^i)|^2 \\
     \label{eq: loss_Stefan_1D1P_L_s0}
    & \mathcal{L}_{s_0}(\bm{\beta})   =  \frac{1}{N} \sum_{i=1}^{N} |s_{\bm \beta}(0) - s_0|^2.
\end{align}
Here $N$ is the batch size and $\{x^i, t^i\}_{i=1}^N$ is a set of collocation points that are randomly sampled inside the computation domain $\Omega^*=[0,1] \times [0, 1] \supset \Omega$ in each training iteration. Since $s(t)$ is unknown in our problem setup, we cannot simply sample points in $\Omega = (0, s(t)) \times (0, T)$. One way to deal with this technical issue is to
sample collocation points in the whole computation domain $\Omega^*$, which is a-priori known, but then just focus on physical solutions by restricting the predicted solution $u_{\bm \theta}(x,t)$ to the domain $(0, s_{\bm \beta}(t)) \times (0,T)$.  

\begin{figure}
    \centering
    \includegraphics[width=0.8\textwidth]{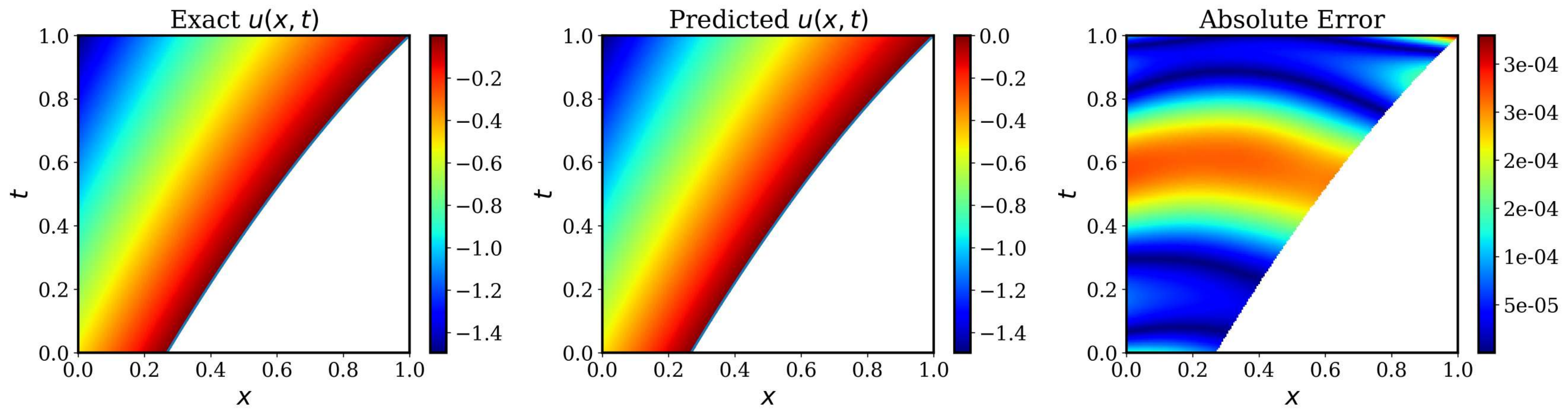}
    \caption{{\em Direct one-dimensional one-phase  Stefan problem:} Exact temperature solution $u(x,t)$ versus the predicted solution by training a  physics-informed neural network using the hyper-parameters summarized in table \ref{tab: hyper-parameters}. The relative $L^2$-error is $2.07e-04$.}
    \label{fig:Stefan_1D1P_direct_pred_u}
\end{figure}

\begin{figure}
    \centering
    \includegraphics[width=0.5\textwidth]{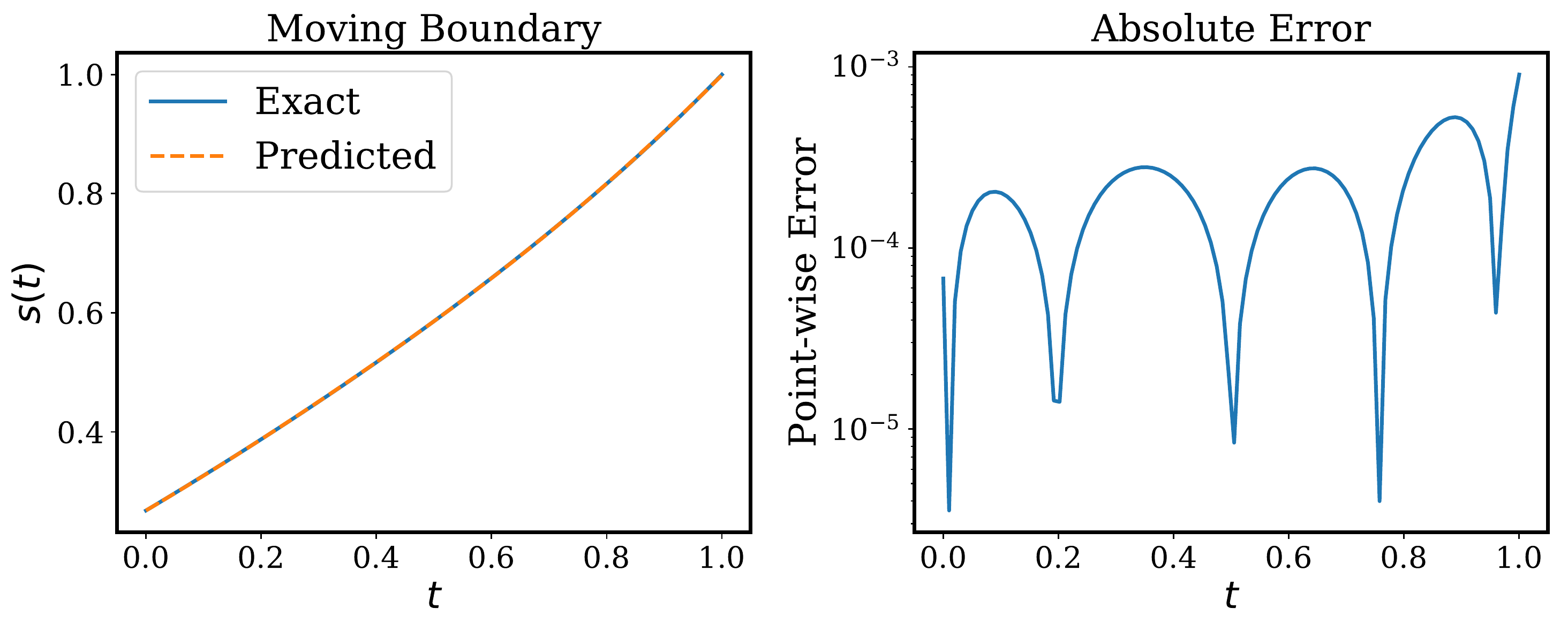}
    \caption{{\em Direct one-dimensional one-phase Stefan problem:} {\em Left:} Comparison of the predicted and exact moving boundary $s(t)$ by training a  physics-informed neural network using the hyper-parameters summarized in table \ref{tab: hyper-parameters}. {\em Right:} The absolute point-wise error between the exact free boundary and the predicted free boundary for $t \in [0,1]$. The relative $L^2$ error is $4.19e-04$.}
    \label{fig:Stefan_1D1P_direct_pred_mb}
\end{figure}

Figures \ref{fig:Stefan_1D1P_direct_pred_u}
and \ref{fig:Stefan_1D1P_direct_pred_mb} present visual comparisons against exact temperature solution $u(x,t)$ and the moving boundary $s(t)$,  as well as point-wise absolute errors.  As it can be seen, 
the approximations are in good agreement with the exact solutions. This figure indicates that our framework is able to obtain accurate predictions of both temperature solution and the moving boundary with a relative $L^2$ error of $2.07e-4$ and $4.19e-04$, respectively.

\subsubsection{Inverse one-dimensional one-phase Stefan problems: Type I}
\label{sec: Stefan_1D1P_inverse_I}

In next two subsections, we will reformulate the last example to introduce two different classes of inverse Stefan problems, as discussed in section \ref{sec: Stefan_1D1P}. The first type of inverse problem we investigate is the so-called moving boundary design problem. Recall that the objective here is to find the solution satisfying equations \ref{eq: Stefan_1D1P_eq} - \ref{eq: Stefan_1D1P_SNc}, as well as recover the Dirichlet and Neumann boundary condition \ref{eq: Stefan_1D1P_ubc&uNc} given the information of free boundary $s(t)$. 
\begin{align}
    \label{eq: Stefan_1D1P_ubc&uNc}
    & u(0,t) = -\frac{1}{2} - t, \quad   \frac{\partial u}{\partial x}(0, t) = 2 , \quad t \in [0, 1]
\end{align}

Note that the free boundary $s(t)$ is known and the only unknown latent variable is the temperature distribution $u(x,t)$. This means that we can deal with this problem by directly leveraging the standard framework of physics-informed neural networks \cite{raissi2019physics}. That is, one can parametrize the unknown solution $u(x,t)$ by a fully-connected neural network $u_{\bm \theta}$ and then train its parameters $\bm{\theta}$ by minimizing the loss function formulated as follows 
\begin{align}
\label{eq: loss_Stefan_1D1P_inverse_I}
 \mathcal{L}({\bm \theta}) =   \mathcal{L}_r({\bm \theta}) +  \mathcal{L}_{u_0}({\bm \theta}) +  \mathcal{L}_{u_{Nc}}({\bm \theta})  +  \mathcal{L}_{s_{bc}}({\bm \theta}) +  \mathcal{L}_{s_{Nc}}({\bm \theta}),
\end{align}
where $\mathcal{L}_r(\theta)$, $\mathcal{L}_{u_0}(\theta)$ and $\mathcal{L}_{u_{Nc}}(\theta)$ are exactly the same as in equations \ref{eq: loss_Stefan_1D1P_L_r} - \ref{eq: loss_Stefan_1D1P_L_uNc}, and
\begin{align}
     & \mathcal{L}_{s_{bc}}({\bm \theta}) = \frac{1}{N} \sum_{i=1}^{N} |u_{\bm \theta}(s(t^i), t^i)|^2 \\ 
     &  \mathcal{L}_{s_{Nc}}({\bm \theta}) = \frac{1}{N} \sum_{i=1}^{N} |\frac{\partial u_{\bm \theta}}{\partial x}(s(t^i), t^i) - h(t^i)|^2.
\end{align}

\begin{figure}
     \centering
     \begin{subfigure}[b]{0.8\textwidth}
         \centering
         \includegraphics[width=\textwidth]{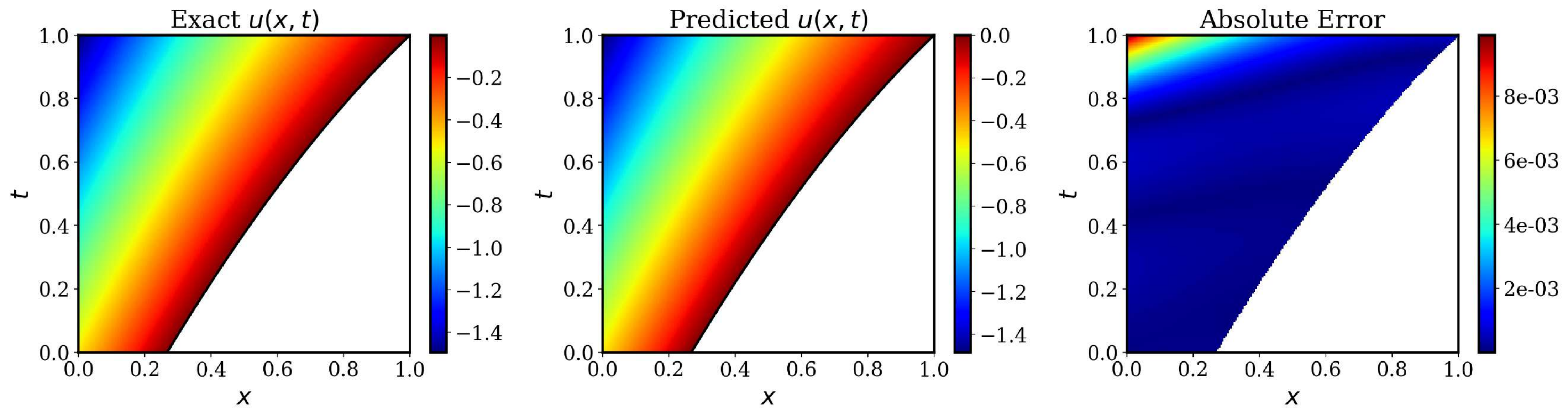}
         \caption{}
         \label{fig: Stefan_1D1P_inverse_I_pred_u}
     \end{subfigure}
     \begin{subfigure}[b]{0.5\textwidth}
         \centering
         \includegraphics[width=\textwidth]{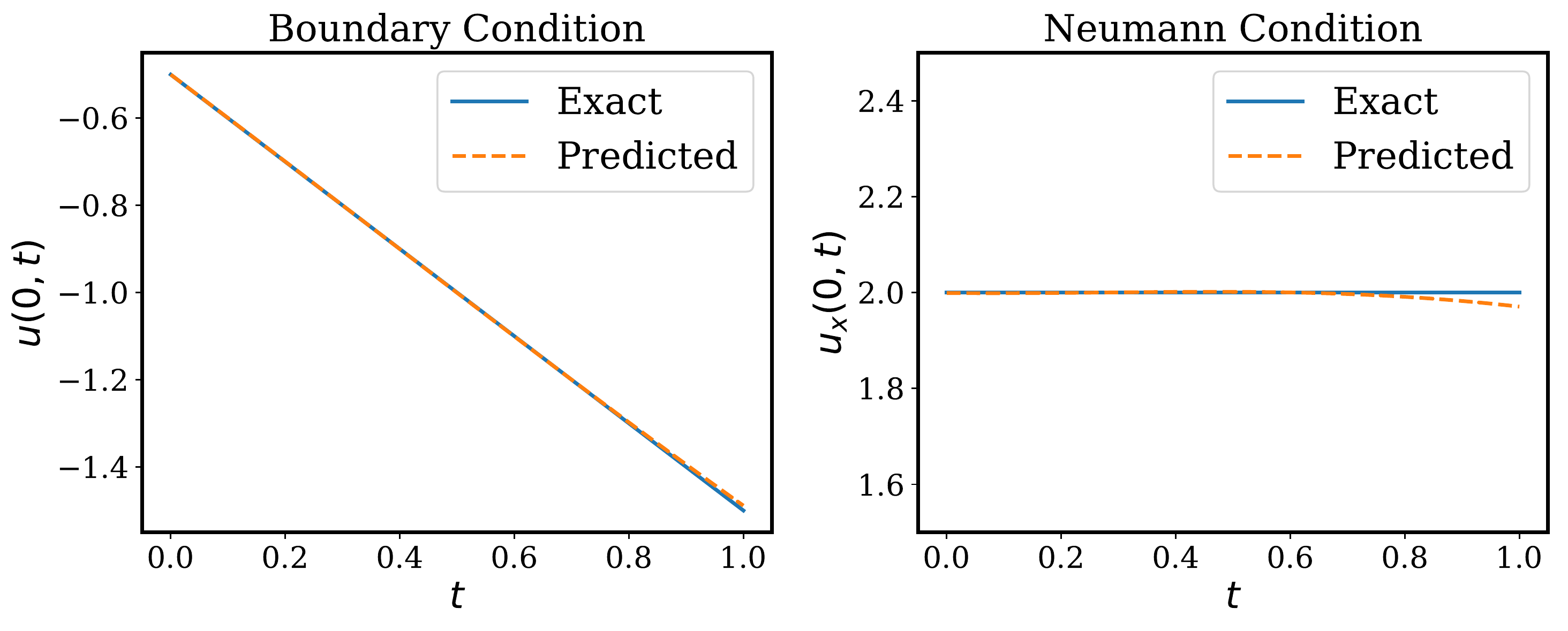}
         \caption{}
         \label{fig: Stefan_1D1P_inverse_I_pred_bc}
     \end{subfigure}
    \caption{{\em Inverse one-dimensional one-phase  Stefan problem Type I: } (a)  Exact temperature solution $u(x,t)$ versus the predicted solution by training a conventional physics-informed neural network model using the hyper-parameters listed in table \ref{tab: hyper-parameters} (relative $L^2$-error: 2.59e-03). (b) {\em Left:} The comparison between the exact $u(0,t)$ and the predicted $u_\theta(0,t)$ (relative $L^2$-error: 3.37e-03);   {\em Right:} The comparison between the exact $u_x(0,t)$ and the predicted $\frac{\partial u_\theta}{\partial x}(0,t)$ (relative $L^2$-error: 5.25e-03).}
    \label{fig: Stefan_1D1P_inverse_I}
\end{figure}

A visual assessment of the obtained predictions against the exact solution $u(x,t)$ is shown in figure \ref{fig: Stefan_1D1P_inverse_I_pred_u}. One can observe that the predicted temperature distribution agrees with the exact one very well and the relative $L^2$ error is $2.59e-03$.  In particular, we present a comparison between the exact and the predicted boundary and Neumann conditions in figure \ref{fig: Stefan_1D1P_inverse_I_pred_bc} where $\frac{\partial u_{\bm \theta}}{\partial x}$ can be computed via automatic differentiation \cite{baydin2017automatic}. From the figure, although the derivative $u_x(0,t)$ becomes less accurate as $t$ goes to final time $T=1$, we see that both $u(0,t)$ and 
$u_x(0,t)$ are recovered relatively well with relative $L^2$ error $3.37e-03$ and $5.25e-3$, respectively. This indicates that the proposed framework is capable of accurately tackling the traditional moving boundary design problem. 

\subsubsection{Inverse one-dimensional one-phase Stefan problems: Type II}
\label{sec: Stefan_1D1P_inverse_II}

In this subsection, we consider a second type of inverse Stefan problems, which we first discussed in section \ref{sec: Stefan_1D1P} and summarized in table \ref{tab: Stefan_1D1P}. Specifically, given a small number of measurements  $\{u(x^j_{\text{data}}, t^j_{\text{data}})\}_{j=1}^M$, we want to infer the solution $u(x,t)$ that satisfies the heat transfer equation \ref{eq: Stefan_1D1P_eq} inside the  domain $\Omega$, as well as identify the unknown free boundary $s(t)$ satisfying initial and boundary conditions \ref{eq: Stefan_1D1P_Sbc}, \ref{eq: Stefan_1D1P_SNc}, and \ref{eq: Stefan_1D1P_Sic}. Here we emphasize that we do not require any initial or boundary conditions, which is close to many realistic applications where such information may be hard to obtain.

To solve this problem, we still follow the methodology outlined in section \ref{sec: Stefan_1D1P_direct}. More specifically, 
the temperature $u(x, t)$ and the moving boundary $s(t)$ are approximated by two independent fully-connected neural networks, $u_{\bm \theta}(x, t) $ and $s_{\bm \beta}(t)$, respectively. Then one can obtain the predicted solutions by minimizing a physics-informed loss of the form
\begin{align}
    \label{eq: loss_Stefan_1D1P_inverse_II}
    \mathcal{L}({\bm \theta}, {\bm \beta}) =  \mathcal{L}_{\text{data}}({\bm \theta}) +  \mathcal{L}_r({\bm \theta}) +  \mathcal{L}_{s_{bc}}({\bm \theta}, {\bm \beta}) +  \mathcal{L}_{s_{Nc}}({\bm \theta}, {\bm \beta})  +   \mathcal{L}_{s_0}({\bm \beta}),
\end{align}
where $\mathcal{L}_r({\bm \theta}), \mathcal{L}_{s_{bc}}({\bm \theta}, {\bm \beta}), \mathcal{L}_{s_{Nc}}({\bm \theta}, {\bm \beta})$ and 
$\mathcal{L}_{s_0}({\bm \beta})$ are the exactly same as in equations \ref{eq: loss_Stefan_1D1P_L_r}, \ref{eq: loss_Stefan_1D1P_L_sbc} - \ref{eq: loss_Stefan_1D1P_L_s0}, and
\begin{align}
    \mathcal{L}_{\text{data}} = \frac{1}{M} \sum_{j=1}^M |u_\theta(x^j_{\text{data}}, t^j_{\text{data}}) - u(x^j_{\text{data}}, t^j_{\text{data}})|^2.
\end{align}

\begin{figure}
     \centering
     \begin{subfigure}[b]{0.8\textwidth}
         \centering
         \includegraphics[width=\textwidth]{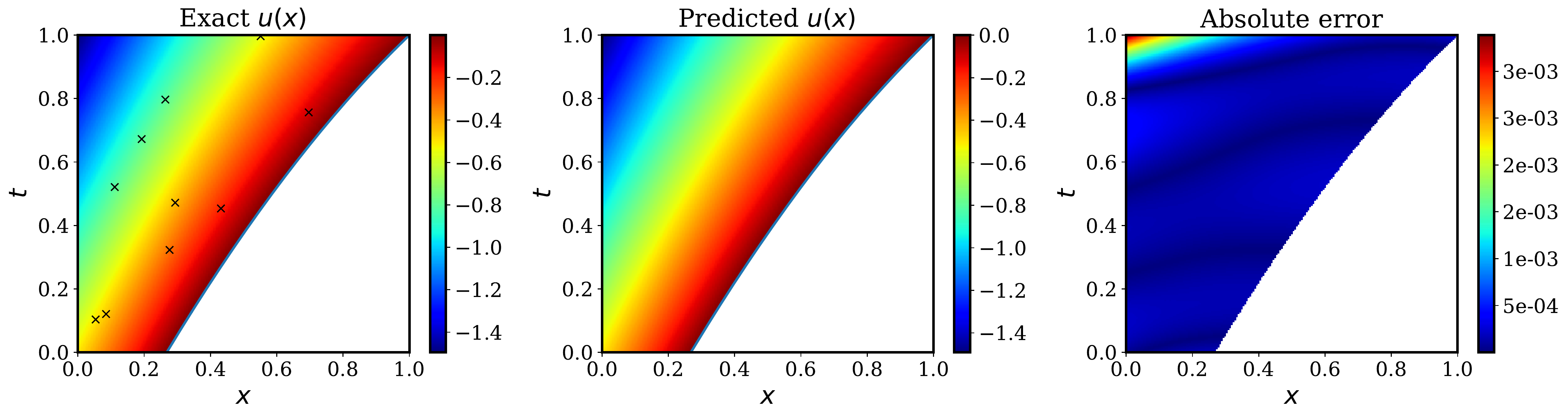}
         \caption{}
         \label{fig: Stefan_1D1P_inverse_II_pred_u}
     \end{subfigure}
     \begin{subfigure}[b]{0.6\textwidth}
         \centering
         \includegraphics[width=\textwidth]{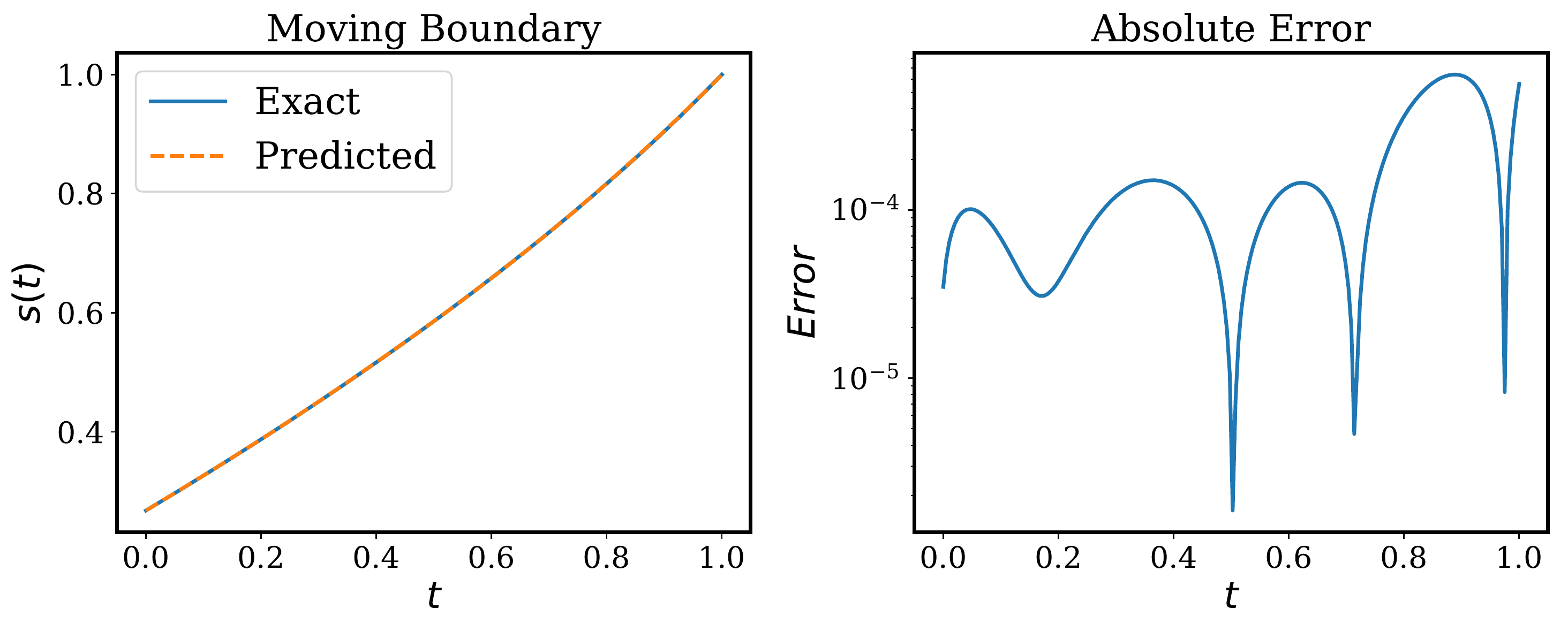}
         \caption{}
         \label{fig: Stefan_1D1P_inverse_II_pred_mb}
     \end{subfigure}
        \caption{{\em Inverse one-dimensional one-phase  Stefan problem Type II: } (a)  Exact temperature solution  versus the predicted solution. The relative $L^2$ error: 5.46e-04. (b) {\em Left:} Comparison of the and exact predicted moving boundary. {\em Right:} The absolute error between the exact  and the predicted free boundary for $t\in [0, 1]$. The $L^2$ error is $3.81e-04$. All plots are obtained using physics-informed neural networks with the hyper-parameters summarized in table \ref{tab: hyper-parameters}.}
        \label{fig:Stefan_1D1P_inverse_II}
\end{figure}

First, we train our model with $M=10$ sparse noise-free measurements of the temperature, which are randomly sampled in $\Omega = (0, s(t)) \times [0, 1]$. A visual assessment of the measurement locations, the resulting predictions for the solution $u(x,t)$ and the moving boundary $s(t)$ is displayed in figure \ref{fig:Stefan_1D1P_inverse_II}.
It can be observed that the proposed methods produce very good predictions for the temperature distributions as well as the moving boundary. 

\begin{figure}
     \centering
     \begin{subfigure}[b]{0.35\textwidth}
         \centering
         \includegraphics[width=\textwidth]{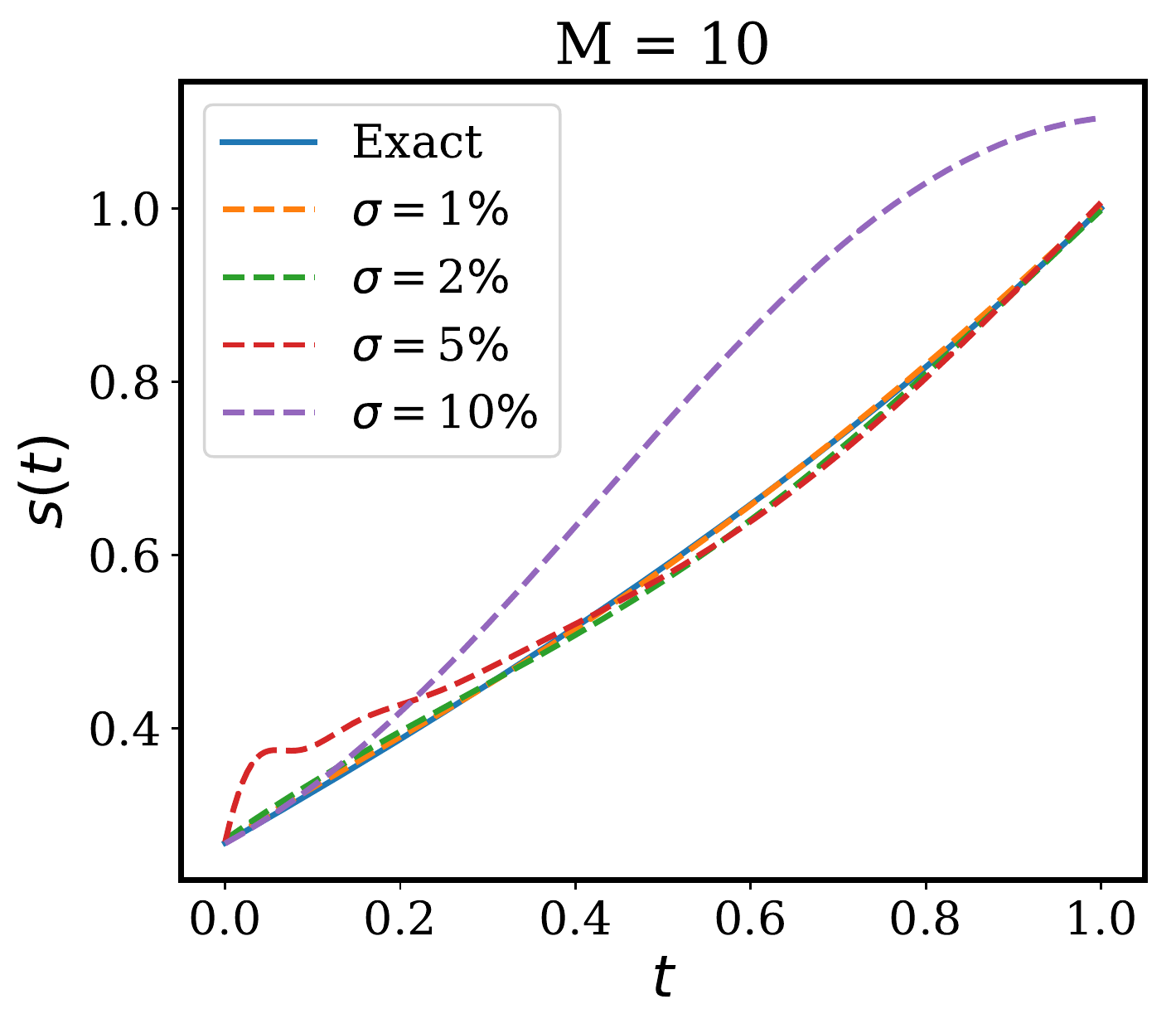}
         \caption{}
         \label{fig: Stefan_1D1P_inverse_II_pred_mb_M_10}
     \end{subfigure}
     \begin{subfigure}[b]{0.35\textwidth}
         \centering
         \includegraphics[width=\textwidth]{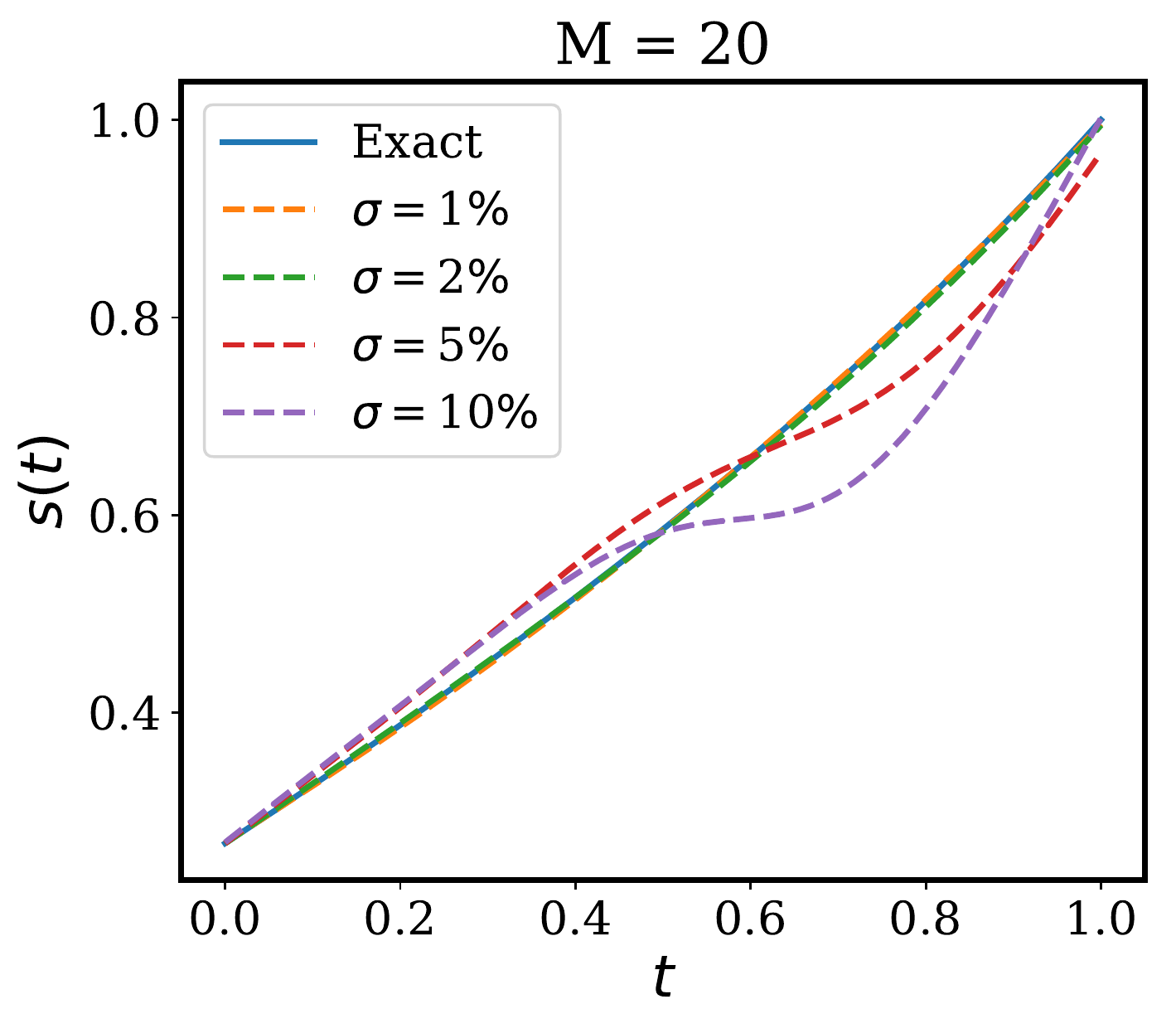}
         \caption{}
         \label{fig:Stefan_1D1P_inverse_II_pred_mb_M_20}
     \end{subfigure}
     \begin{subfigure}[b]{0.35\textwidth}
         \centering
         \includegraphics[width=\textwidth]{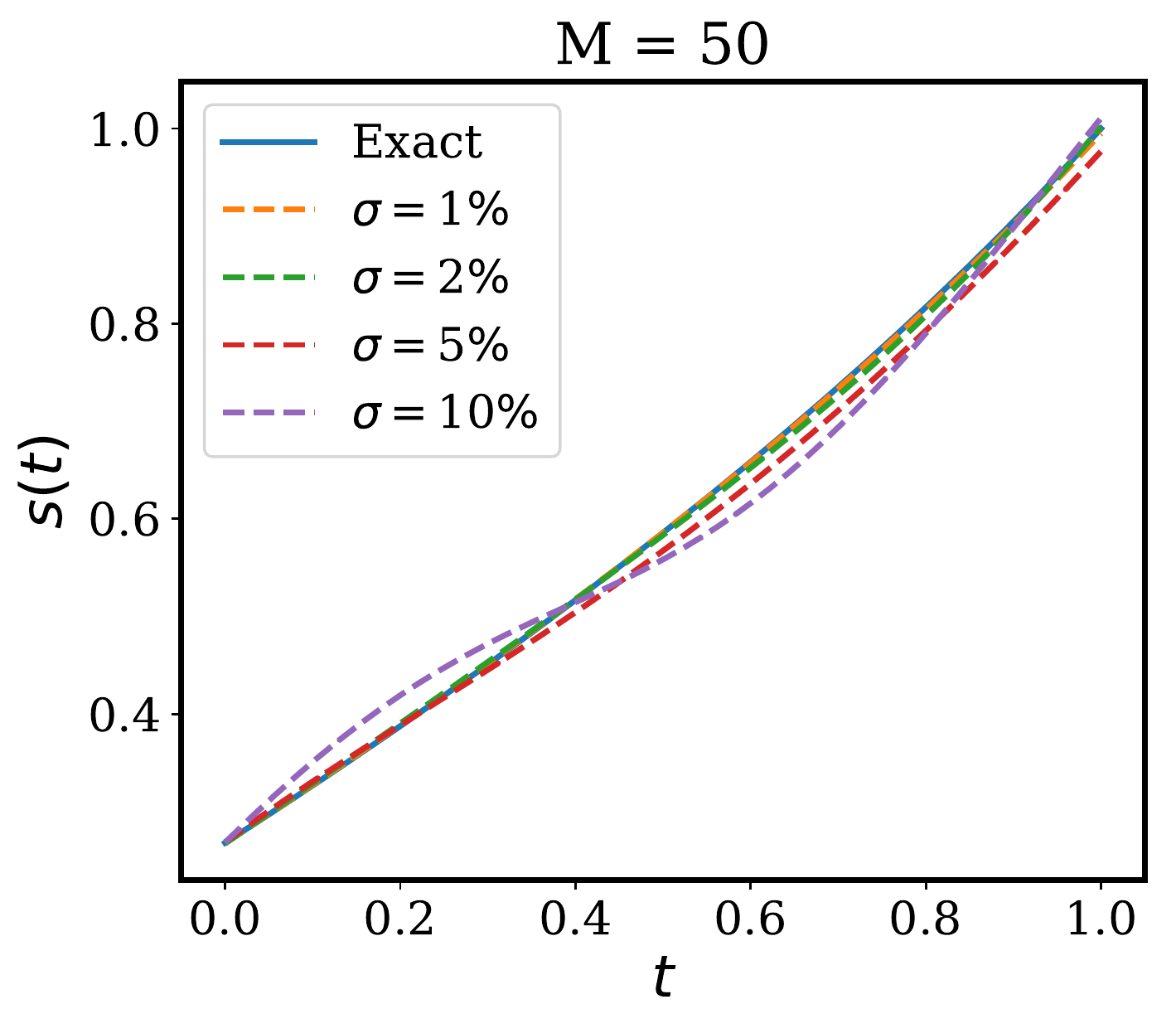}
         \caption{}
         \label{fig: Stefan_1D1P_inverse_II_pred_mb_M_50}
     \end{subfigure}
          \begin{subfigure}[b]{0.35\textwidth}
         \centering
         \includegraphics[width=\textwidth]{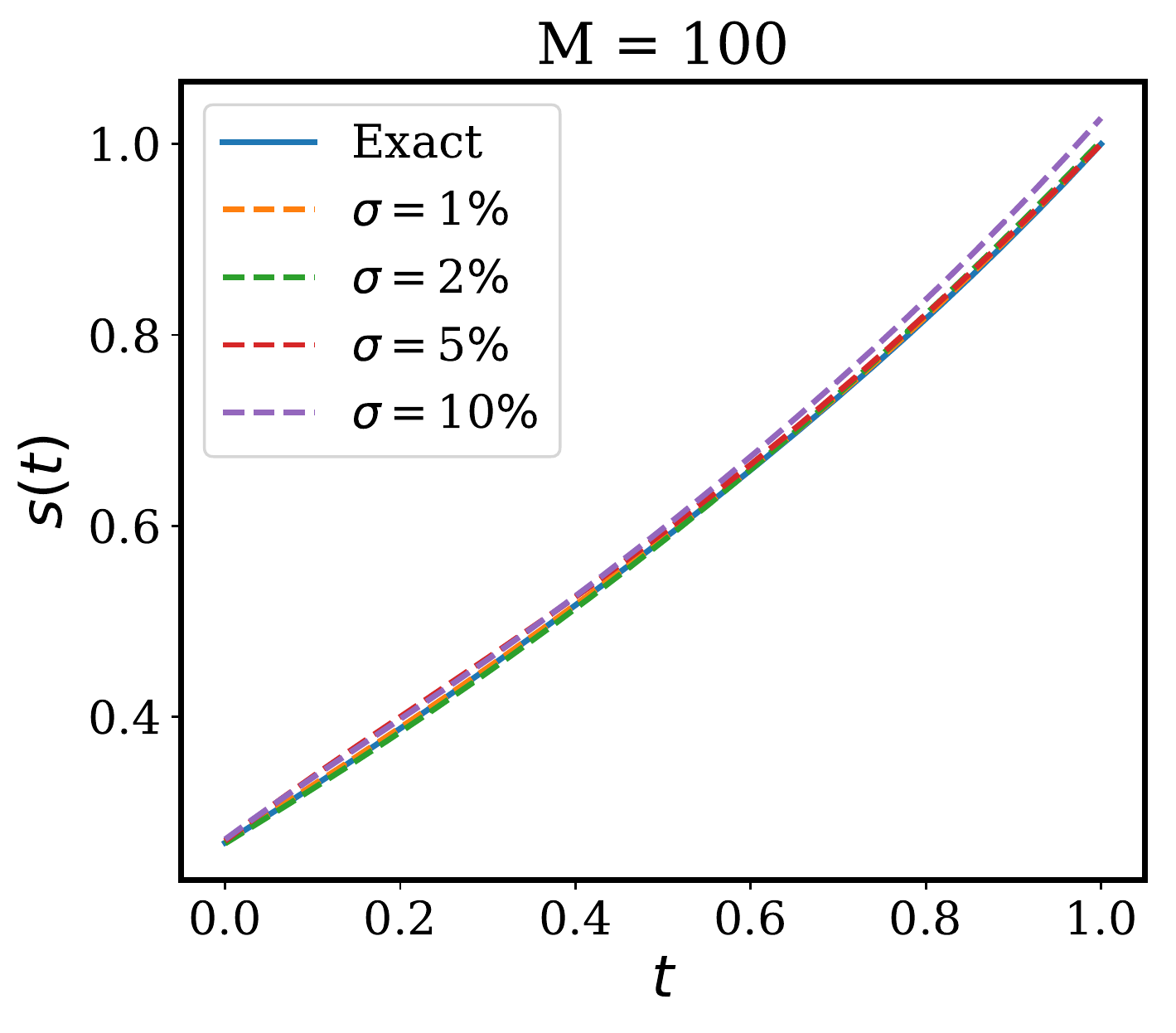}
         \caption{}
         \label{fig: Stefan_1D1P_inverse_II_pred_mb_M_100}
     \end{subfigure}
        \caption{{\em Inverse one-dimensional one-phase inverse Stefan problem Type II: } Comparison between the exact and the predicted moving boundary corresponding to different number of data points $M \in \{10, 20, 50, 100\}$ and with different noise levels $\sigma \in \{1\%, 2\%, 5\%, 10\%\}$. All plots are obtained using physics-informed neural networks with hyper-parameters summarized in table \ref{tab: hyper-parameters}.
        }
        \label{fig:Stefan_1D1P_inverse_II_pred_mb_different_M}
\end{figure}

In the next step, we aim to study the sensitivity of the proposed framework to the number of data  $M$, as well as the variance $\sigma$ of the uncorrelated Gaussian noise that may corrupts the temperature measurements, taken to be proportional to the maximum value attained by the solution, i.e.,
\begin{align}
\label{eq: noise}
    \sigma = \delta \times \|u(x,t)\|_{L^\infty(\Omega)}
\end{align}
where $\delta$ is the level of noise. 
To this end, we systematically analyze the performance of the proposed methods and quantify their predictive accuracy for different number of data points $M$ and different noise levels. The results are summarized in table \ref{tab:Stefan_1D1P_inverse_II_pred_u} and table \ref{tab:Stefan_1D1P_inverse_II_pred_mb}. Also, a detailed visual assessment of predictive moving boundary is shown in \ref{fig:Stefan_1D1P_inverse_II_pred_mb_different_M}.  One can conclude that the proposed framework remains robust with respect to the noise level in the data, and yields a reasonable identification accuracy, even for noise corruptions up to $10 \%$. Our experience so far indicates that the negative consequences of more noise in the data can be remedied to some extent by obtaining more data. In addition, in the noise-free case, only a handful of measurements is required to attain good predictive accuracy for both the  temperature solution and the moving boundary. 

\begin{table}[]
\renewcommand*{\arraystretch}{1.6}
    \centering
    \begin{tabular}{|c|c|c|c|c|c|}
\hline
\diagbox{Data Points M}{Noise $\sigma$}  & $\sigma =0\%$      & $\sigma =1\%$ &  $\sigma =2\%$ &  $\sigma =5\%$&  $\sigma =10\%$ \\ \hline
M = 10 &8.98e-04 & 6.57e-03   & 2.89e-02    &  7.34e-02   &  4.24e-01    \\ \hline
M = 20 &4.72e-04 & 1.06e-02    & 1.10e-02    &  1.05e-01    & 2.01e-01     \\ \hline
M = 50 &4.55e-04 & 7.62e-03     &  2.48e-02   & 5.15e-02    & 8.91e-02     \\ \hline
M = 100&4.04e-04 & 3.73e-03    & 1.05e-02    & 2.02e-02    & 4.87e-02     \\ \hline
\end{tabular}
    \caption{{\em Inverse one-dimensional one-phase Stefan problem Type II:} Relative $L^2$ errors of predicted solution $u(x,t)$ with different number of data points and different noise level $\sigma$.}
    \label{tab:Stefan_1D1P_inverse_II_pred_u}
\end{table}

\begin{table}[]
\renewcommand*{\arraystretch}{1.6}
    \centering
    \begin{tabular}{|c|c|c|c|c|c|}
\hline
\diagbox{Data Points M}{Noise $\sigma$}  & $\sigma =0\%$      & $\sigma =1\%$ &  $\sigma =2\%$ &  $\sigma =5\%$&  $\sigma =10\%$ \\ \hline
M = 10 &3.12e-04 & 4.08e-03  & 1.69e-02    &  4.53e-02   &  2.30e-01  \\ \hline
M = 20 &3.47e-04 & 2.67e-03    & 5.96e-03    & 5.39e-02    &9.11e-02     \\ \hline
M = 50 &4.37e-04 &2.76e-03     &  7.64e-03   &2.72e-02   & 4.17e-02     \\ \hline
M = 100&3.10e-04 & 2.07e-03   & 5.48e-03   & 1.22e-02    & 2.40e-02     \\ \hline
\end{tabular}
    \caption{{\em Inverse one-dimensional two-phase Stefan problem Type II:} Relative $L^2$ errors of predicted free boundary $s(t)$ with different number of data points and different noise level $\sigma$.}
    \label{tab:Stefan_1D1P_inverse_II_pred_mb}
\end{table}

\subsection{One-dimensional two-phase Stefan problems}
\label{sec: Stefan_1D2P}
Our second case study is focused on extending the proposed methods to one-dimensional, two-phase Stefan problems. In contrast to the classical one-phase case, the numerical solution of direct and inverse two-phase Stefan problems has received very little attention in the literature, see \cite{slota2008solving,magenes1989theoretical,nochetto1991adaptive, johansson2013meshless}. Moreover, in most existing references, two-phase inverse Stefan problems are solved assuming that the moving boundary is known, as the presence of unknown boundaries renders the numerical solution of two-phase Stefan problems much more difficult. Motivated by this capability gap, in this section we focus on applying the proposed deep learning framework to direct and inverse one-dimensional two-phase Stefan problems under the assumption that the moving boundary is unknown.

Let us start by giving the mathematical formulation of a classical one-dimensional two phase Stefan problem \cite{johansson2013meshless} as follows. Assume a moving boundary $s(t) \subset [0,L]$ for $t \in [0,T]$ and consider a two-phase rectangular domain, 
\begin{align}
    \Omega = \{(x,t) : (0, L) \times (0, T]\}
\end{align}
which is subdivided by a moving interface into two sub domains
\begin{align}
    & \Omega_1 = \{ (x,t) \in \Omega : 0< x < s(t), t \in (0,T] \} \\
    & \Omega_2 = \{ (x,t) \in \Omega : s(t)< x < L, t \in (0,T] \}.
\end{align}
Then, the latent temperature distributions within each of the two phases, $u_1(x,t)$ and $u_2(x,t)$, respectively, satisfy a heat equation
\begin{align}
     \label{eq: Stefan_1D2P_eq}
    &\frac{ \partial u_i }{\partial t} = k_i \frac{\partial^2 u_i}{\partial x^2}, \quad (x, t) \in \Omega_i, \quad i = 1,2,
\end{align}
subject to the initial and boundary conditions 
\begin{align}
     \label{eq: Stefan_1D2P_ic}
    & u_i(x, 0) =u_0^{(i)}(x,0 ), \quad x \in [0, s(0)], \quad i=1,2 \\
    \label{eq: Stefan_1D2P_bc}
    & u_1(0, t) =  g_1(t)  , \quad t \in [0, 1] \\
    & u_2(l, t) =  g_2(t), \quad t \in [0, 1].
\end{align}
Also, the initial and boundary Stefan conditions are given by
\begin{align}
    \label{eq: Stefan_1D12P_Sbc}
    & u_1(s(t), t) = u_2(s(t), t) = u^*, \quad t \in [0, 1] \\
    \label{eq: Stefan_1D12P_SNc} 
    & s'(t) = \alpha_1 \frac{\partial u_1}{\partial x}(s(t), t) + \alpha_2 \frac{\partial u_2}{\partial x}(s(t), t),  \quad t \in [0, 1] \\
    \label{eq: Stefan_1D12P_Sic}
    &s(0) = s_0.
\end{align}

For direct one-dimensional two-phase problem, our goal is to infer the unknown  temperature distributions $u_i(x,t)$, as well as the moving interface $s(t)$ that satisfies equations \ref{eq: Stefan_1D2P_eq} - \ref{eq: Stefan_1D12P_Sic}.


To demonstrate the generality of our framework, we also formulate another type of inverse Stefan problem for the two-phase case, which is different from
the moving boundary problem design problem considered in section \ref{sec: Stefan_1D1P_inverse_I}. That is, instead of assuming knowledge of the boundary conditions \ref{eq: Stefan_1D2P_bc1} and \ref{eq: Stefan_1D2P_bc2},
we provide additional data at the final time, i.e., equation \ref{eq: Stefan_1D12P_ft} 
\begin{align}
    \label{eq: Stefan_1D12P_ft}
    u_i(x,T) = h_i(x),  \quad i =1,2.
\end{align}
Then our objective is to find the temperature solutions $u_1(x,t)$ and $u_2(x,t)$, and the unknown free boundary $s(t)$, satisfying \ref{eq: Stefan_1D2P_eq},\ref{eq: Stefan_1D2P_ic}, \ref{eq: Stefan_1D12P_Sbc} - \ref{eq: Stefan_1D12P_Sic} and \ref{eq: Stefan_1D12P_ft}. In the context of one-dimensional two-phase Stefan problems, we will refer to this class of inverse problems as Inverse Type I.

Furthermore, the second type of one-dimensional one-phase Stefan inverse problem in the section \ref{sec: Stefan_1D2P_inverse_II} can be naturally extended to the one-dimensional two-phase case. Similarly, our goal is to infer the triple solution $\{u_1(x,t), u_2(x,t), s(t)\}$ satisfying the equations \ref{eq: Stefan_1D2P_eq}, \ref{eq: Stefan_1D12P_Sbc} - \ref{eq: Stefan_1D12P_Sic}, assuming that some measurements of temperature are given inside the domain $\Omega$.
Table \ref{tab: Stefan_1D2P} provides a summary of all one-dimensional two-phase Stefan problems considered in this section.

 \begin{table}
    \centering
     \renewcommand*{\arraystretch}{1.8}
   \begin{tabular}{|c|c|c|c|}
\hline
                & Equations  & Observed & Latent \\ \hline
Direct          & \ref{eq: Stefan_1D2P_eq} - \ref{eq: Stefan_1D12P_Sic}.       &  $u_0^{(i)}(x), g_i(t), k_i, \alpha_i, u^*, s_0 \quad i=1,2 $ & $u_i(x,t), s(t), \quad i=1,2$       \\ \hline
Inverse Type I  &  \ref{eq: Stefan_1D2P_eq},\ref{eq: Stefan_1D2P_ic}, \ref{eq: Stefan_1D12P_Sbc} - \ref{eq: Stefan_1D12P_Sic}, \ref{eq: Stefan_1D12P_ft}          & $u_0^{(i)}(x), h_i(x), k_i, \alpha_i, u^*, s_0   \quad i=1,2$ & $u_i(x,t), s(t), \quad i=1,2  $     \\ \hline
Inverse Type II &\ref{eq: Stefan_1D2P_eq}, \ref{eq: Stefan_1D12P_Sbc} - \ref{eq: Stefan_1D12P_Sic},       & $\{(x^j_{\text{data}}, t^j_{\text{data}}), u^j \}_{j=1}^M, k_i, \alpha_i, u^*, s_0$    &  $u_i(x,t), s(t), \quad i=1,2  $    \\ \hline
\end{tabular}
    \caption{{\em One-dimensional two-phase Stefan problem:} Summary of conditions and the objective of each type of Stefan problem formulated in section \ref{sec: Stefan_1D2P}.}
    \label{tab: Stefan_1D2P}
\end{table}

\subsubsection{Direct one-dimensional two-phase Stefan problems}
\label{sec: Stefan_1D2P_direct}

Here we consider a concrete example of the direct one-dimensional two-phase Stefan problem \cite{johansson2013meshless}, in a computational domain $ \Omega = [0, 2] \times [0, 1]$ with given coefficients $k_1 = 2$, $k_2=1$, $\alpha_1 =-2$, $\alpha_2 =1$, and $s_0 = 1/2$. Moreover, we assume the following initial and boundary conditions 
\begin{align}
     \label{eq: Stefan_1D2P_ic1}
    & u_0^{(1)}(x,0) =2 [\exp(\frac{1/2 - x}{2}) -1] , \quad x \in [0, s(0)] \\
    \label{eq: Stefan_1D2P_ic2}
    & u_0^{(2)}(x,0 ) = \exp(1/2 - x) - 1 , \quad x \in [s(0), 2] \\
    \label{eq: Stefan_1D2P_bc1}
    &  g_1(t) = 2 [\exp(\frac{t + 1/2}{2}) -1] , \quad t \in [0, 1] \\
     \label{eq: Stefan_1D2P_bc2}
    &  g_2(t) =\exp(t - 3/2) - 1, \quad t \in [0, 1],
\end{align}
corresponding to an exact solution for the temperature distributions and the moving boundary given by
\begin{align}
    \label{eq: Stefan_1D2P_u1_sol}
    &u_1(x, t) = 2 (\exp((t + 1/2 - x) / 2) - 1) \\
     \label{eq: Stefan_1D2P_u2_sol}
    &u_2(x, t) = \exp((t + 1/2 - x) - 1) \\
    \label{eq: Stefan_1D2P_mb_sol}
    & s(t) = t + 1/2.
\end{align}

Recall that, in direct Stefan problems, we seek to find the temperature distribution in $\Omega_1$ and $\Omega_2$ as well as infer the moving boundary $s(t)$, which satisfies equations \ref{eq: Stefan_1D2P_eq} - \ref{eq: Stefan_1D12P_Sic}.

Similar to the one-phase Stefan problems discussed in section \ref{sec: Stefan_1D1P}, the unknown free boundary $s(t)$ can be approximated by a deep fully-connected neural network $s_{\bm \beta}$. 
However, a key difference here is how we approximate the temperature distribution $u_1(x,t)$ and $u_2(x,t)$ by neural networks. If one employs two independent neural networks with a single output to approximate the temperature distributions $u_1(x,t)$ and $u_2(x,t)$ together in the whole domain $\Omega$, then a difficulty arises in how to compute the normal derivative $\frac{\partial u_1}{\partial x}$ and $\frac{\partial u_2}{\partial x}$ at the interface, leading to difficulties in how to then 
enforce the Stefan Neumann condition \ref{eq: Stefan_1D12P_SNc}. To avoid this problem, we will approximate $u_1(x,t)$ and $u_2(x,t)$  by a single fully-connected neural network with two outputs $u^{(1)}_{\bm \theta}$ and $u^{(2)}_{\bm \theta}(x,t)$, i.e,

\[
    [x, t] \xmapsto[]{u_{\bm \theta}} [u^{(1)}_{\bm \theta}(x,t), u^{(2)}_{\bm \theta}(x,t)]
\]

In this way, the final predicted solution in the whole domain $\Omega$ can be given by
\begin{align}
    \label{eq: Stefan_1D2P_u_pred}
    u_{\bm{\theta}}(x,t) = u^{(1)}_{\bm \theta}(x,t)\textbf{1}_{(0, s_{\bm \beta} (t))}(x) +  u^{(2)}_{\bm \theta}(x,t)\textbf{1}_{(s_{\bm \beta} (t), 1)}(x),
\end{align}
where $\textbf{1}_{(0, s_{\bm \beta} (t))}(x)$ denotes an indicator function taking a value of one if $x\in[0, s_{\bm \beta} (t)]$ and zero otherwise.
Then, the exact temperature solution in the whole domain $\Omega$ can be expressed by
\begin{align}
     \label{eq: Stefan_1D2P_u_exact}
    u_{\text{exact}}(x,t) = u_1(x,t)\textbf{1}_{(0, s (t))}(x) +  u_2(x,t)\textbf{1}_{(s(t), 1)}(x)
\end{align}
The shared parameters $\bm{\theta,\beta}$ can
be learned by minimizing the mean squared error loss
\begin{align}
    \label{eq: loss_1D2P_direct}
    \mathcal{L}(\bm{\theta,\beta}) =  \sum_{k=1}^2         \Big[  \mathcal{L}_r^{(k)}(\bm{\theta}) + \mathcal{L}_{u_0}^{(k)}(\bm{\theta}) +  \mathcal{L}_{u_{bc}}^{(k)}(\bm{\theta})  +  \mathcal{L}_{s_{bc}}^{(k)}(\bm{\theta, \beta}) \Big] +  \mathcal{L}_{s_{Nc}}(\bm{\theta,\beta}) + \mathcal{L}_{s_0}(\bm{\beta}),
\end{align}
where,
\begin{align}
     \label{eq: loss_Stefan_1D2P_L_r}
    &\mathcal{L}_r^{(k)}(\bm{\theta}) = \frac{1}{N} \sum_{i=1}^{N}|\frac{\partial u_{\bm \theta}^{(k)}}{\partial t}(x^i, t^i) -k_i \frac{\partial^2 u_{\bm \theta}^{(k)}}{\partial x^2}(x^i, t^i)|^2    \\
     \label{eq: loss_Stefan_1D2P_L_u0}
    &\mathcal{L}_{u_0}^{(k)}(\bm{\theta}) = \frac{1}{N} \sum_{i=1}^{N} |u_{\bm \theta}^{(k)}(x^i,0) - u_0^{(k)}(x^i)|^2 \\
     \label{eq: loss_Stefan_1D2P_L_ubc}
    & \mathcal{L}_{u_{bc}}^{(k)}(\bm{\theta}) = \frac{1}{N} \sum_{i=1}^N|\frac{\partial u_{\bm \theta}^{(k)}}{\partial x}(0, t^i) - g_k(t^i)|^2  
     \\
     \label{eq: loss_Stefan_1D2P_L_sbc}
    & \mathcal{L}_{s_{bc}}^{(k)}(\bm{\theta,\beta}) = \frac{1}{N} \sum_{i=1}^{N} |u_{\bm \theta}^{(k)}(s_{\bm \beta}(t^i), t^i)|^2,\\
\end{align}
for $k=1,2$, and
\begin{align}
      \label{eq: loss_Stefan_1D2P_L_sNc}
    &  \mathcal{L}_{s_{Nc}}(\bm{\theta,\beta}) = \frac{1}{N} \sum_{i=1}^{N} |2 \frac{\partial u_{\bm \theta}^{(1)}}{\partial x}(s_{\bm \beta}(t^i), t^i) - \frac{\partial u_{\bm \theta}^{(2)}}{\partial x}(s_{\bm \beta}(t^i), t^i) + \frac{d s_{\bm \beta}}{d t}(t^i)|^2\\
     \label{eq: loss_Stefan_1D2P_L_s0}
    & \mathcal{L}_{s_0}({\bm \beta})  =  \frac{1}{N} \sum_{i=1}^{N} |s_{\bm \beta}(0) - s_0|^2.
\end{align}

Again, $N$ denotes the batch size and $(x^i, t^i)$ are collocation points that are randomly sampled inside the entire domain $\Omega$.

\begin{figure}
     \centering
     \begin{subfigure}[b]{0.8\textwidth}
         \centering
         \includegraphics[width=\textwidth]{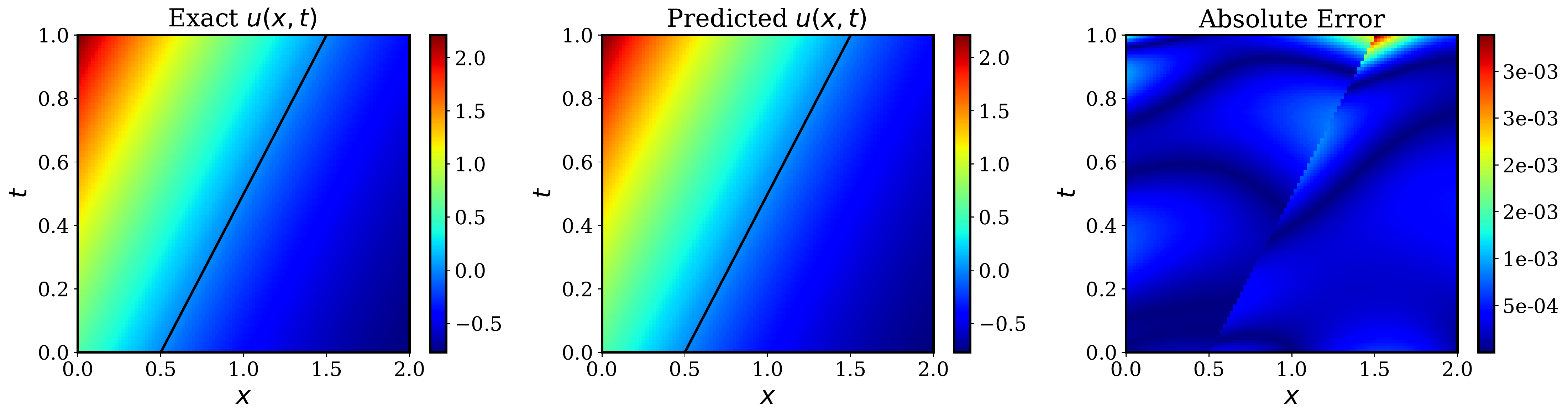}
         \caption{}
         \label{fig:y Stefan_1D2P_direct_pred_u}
     \end{subfigure}
     \begin{subfigure}[b]{0.5\textwidth}
         \centering
         \includegraphics[width=\textwidth]{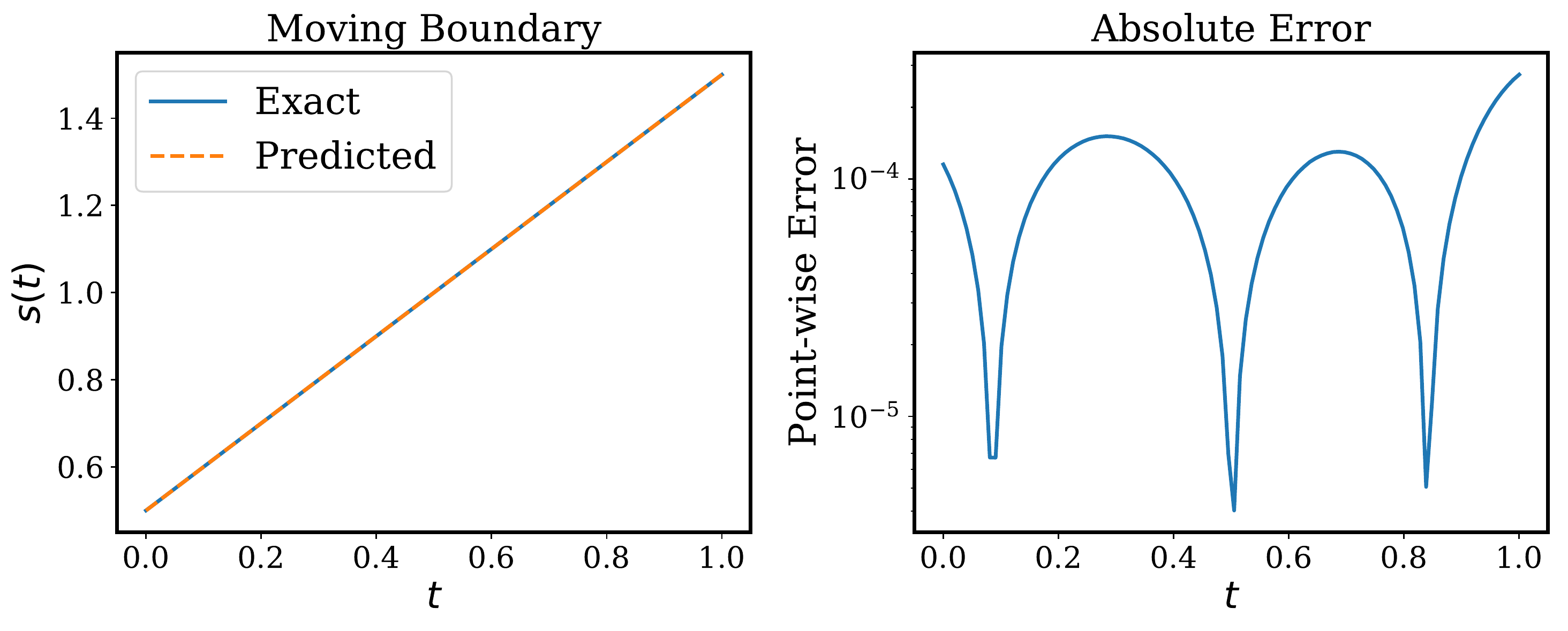}
         \caption{}
         \label{fig: Stefan_1D2P_direct_pred_mb}
     \end{subfigure}
        \caption{{\em Direct one-dimensional two-phase Stefan problem: } (a) Exact temperature solution  versus the predicted solution. The relative $L^2$ error: 5.83e-04. (b) {\em Left:} Comparison of the and exact predicted moving boundary. {\em Right:} The absolute error between the exact  and the predicted free boundary for $t\in [0, 1]$. The $L^2$ error is 2.81e-04. All plots are obtained using physics-informed neural networks with hyper-parameters summarized in table \ref{tab: hyper-parameters}.}
        \label{fig: Stefan_1D2P_direct}
\end{figure}

The results for direct one-dimensional two-phase Stefan problem are depicted in figure \ref{fig: Stefan_1D2P_direct}. In the top
panel, we present a comparison of the exact and the predicted temperature distribution in $\Omega$ along with the exact and predicted moving boundary denoted by the black line. A more detailed comparison of the moving boundary is given in the bottom panel. We observe that both $u_1(x,t)$ and $u_2(x,t)$ have good agreement with the exact temperature distributions with relative $L^2$ error $5.83e-04$ in the entire domain. 
Moreover, the predicted moving boundary also attains high accuracy against the exact solution with the $L^\infty$ error of  $O(10^{-4})$. These observations indicate that the proposed framework is able to effectively tackle this class of direct one-dimensional two-phase Stefan problems.

\subsubsection{Inverse one-dimensional two-phase Stefan problems: Type I}
\label{sec: Stefan_1D2P_inverse_I}
In this section, we try to solve the first type inverse one-dimensional two-phase Stefan problem introduced in section \ref{sec: Stefan_1D2P}, namely Inverse Type I, see table \ref{tab: Stefan_1D2P} for details. Here all values of parameters and the exact solutions are same as in section \ref{sec: Stefan_1D2P_direct}, albeit we will substitute the boundary condition \ref{eq: Stefan_1D2P_bc} by 
providing  additional data at the final time, i.e., equation \ref{eq: Stefan_1D12P_ft1} and  \ref{eq: Stefan_1D12P_ft2}, where
\begin{align}
    \label{eq: Stefan_1D12P_ft1}
     &h_1(x) =  2 (\exp((1 + 1/2 - x) / 2) - 1) ,\quad x \in (0, s(T)) \\
    \label{eq: Stefan_1D12P_ft2}
    & h_2(x) =  \exp((1 + 1/2 - x) - 1,  \quad x \in (s(T), 2),
\end{align}
where we set the final time $T=1$. The the exact solutions for all $\{(u_1(x,t),u_2(x,t), s(t))\}$ are still same as given by equations \ref{eq: Stefan_1D2P_u1_sol} - \ref{eq: Stefan_1D2P_mb_sol}.

We use a single fully-connected neural network with two outputs to approximate $u_1(x,t)$ and $u_2(x,t)$, and another fully-connected network to approximate $s(t)$. We use the same formula as \ref{eq: Stefan_1D2P_u_pred} to obtain the predicted solution $u(x,t)$ in the whole domain $\Omega$, while the physics-informed loss function is now given by
\begin{align}
    \label{eq: loss_Stefan_1D2P_inverse_I}
     \mathcal{L}(\bm{\theta, \beta}) =  \sum_{k=1}^2 \Big[ \mathcal{L}_r^{(k)}({\bm \theta}) + \mathcal{L}_{u_0}^{(k)}(\bm{\theta}) +  \mathcal{L}_{u_{T}}^{(k)}(\bm{\theta})  +  \mathcal{L}_{s_{bc}}^{(k)}(\bm{\theta, \beta}) \Big]  +  \mathcal{L}_{s_{Nc}}(\bm{\theta, \beta}) + \mathcal{L}_{s_0}(\bm{\beta}),  \quad \text{for}  k=1,2.
\end{align}
where all summands are exactly the same as in equations \ref{eq: loss_Stefan_1D2P_L_r} and \ref{eq: loss_Stefan_1D2P_L_ubc}  - \ref{eq: loss_Stefan_1D2P_L_s0} described for the direct problem, except for $\mathcal{L}_{u_T}$ which is now defined as
\begin{align}
    \label{eq: loss_Stefan_1D2P_L_uT}
     \mathcal{L}_{u_{T}}^{(k)} = \frac{1}{N} \sum_{i=1}^{N} | u_{\bm \theta}^{(k)}(x^i, T) - u_k(x^i, T) |^2,
\end{align}

\begin{figure}
     \centering
     \begin{subfigure}[b]{0.8\textwidth}
         \centering
         \includegraphics[width=\textwidth]{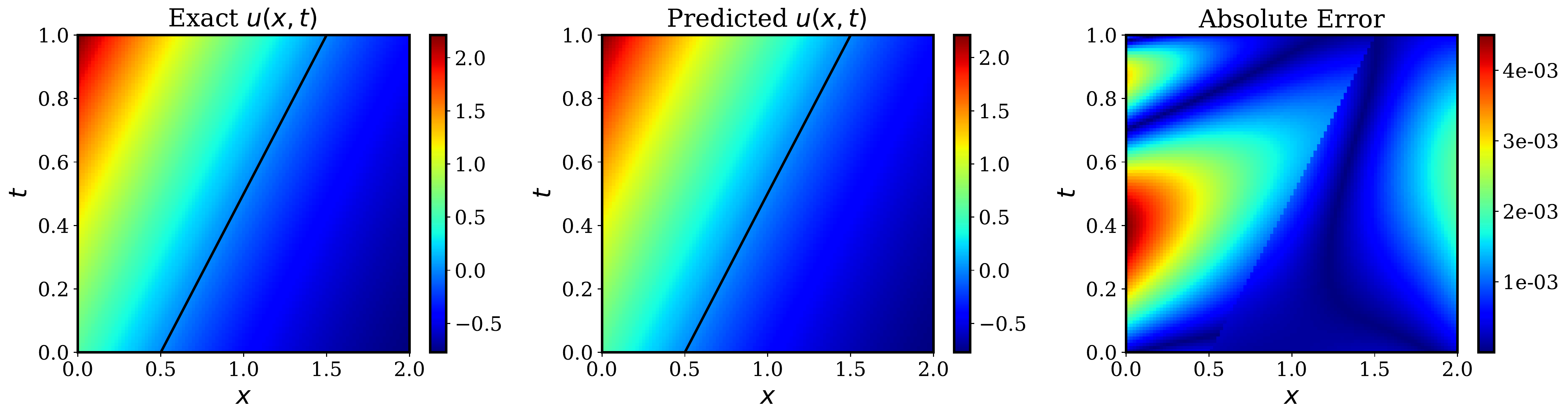}
         \caption{l2 error: 1.91e-03}
         \label{fig: Stefan_1D2P_inverse_I_pred_u}
     \end{subfigure}
     \begin{subfigure}[b]{0.5\textwidth}
         \centering
         \includegraphics[width=\textwidth]{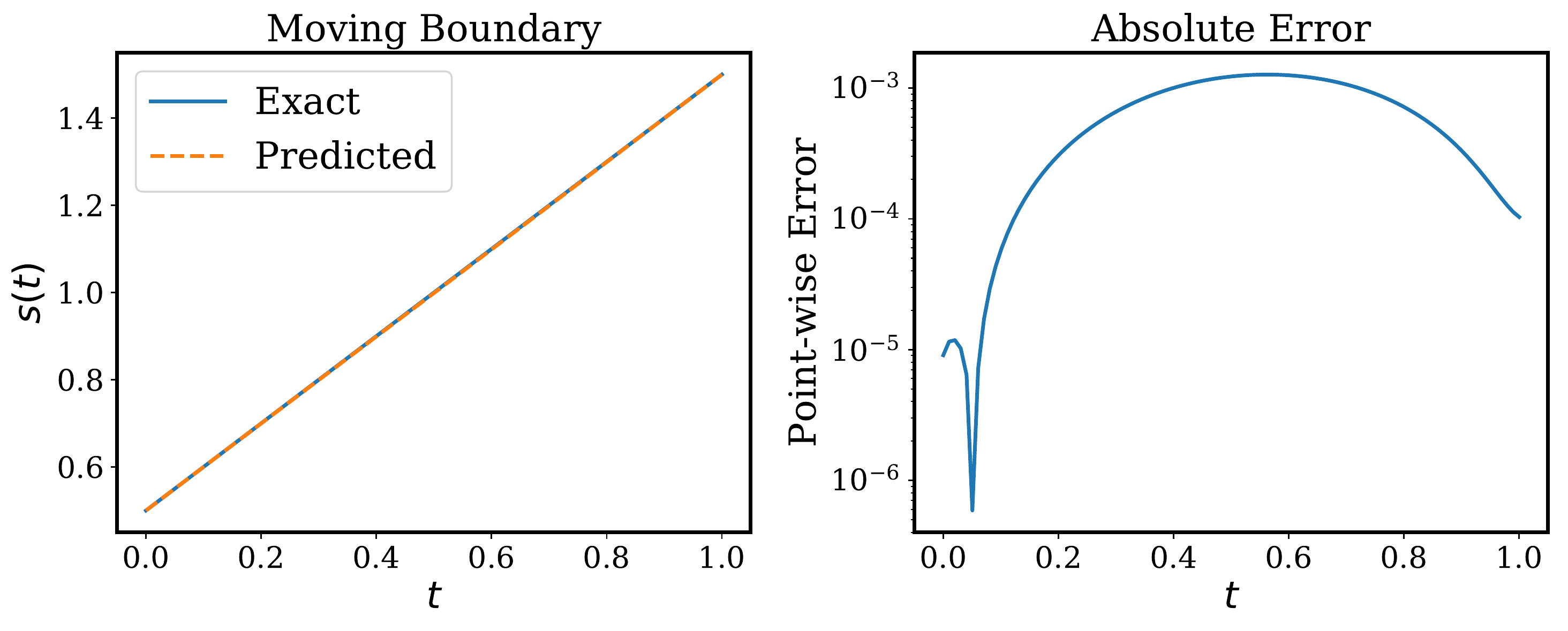}
         \caption{l2 error:7.01e-04}
         \label{fig: Stefan_1D2P_inverse_I_pred_mb}
     \end{subfigure}
        \caption{{\em Inverse one-dimensional two-phase Stefan problem Type I: } (a) Exact temperature solution  versus the predicted solution. The relative $L^2$ error: $1.91e-03$. (b) {\em Left:} Comparison of the and exact predicted moving boundary. {\em Right:} The absolute error between the exact  and the predicted free boundary for $t\in [0, 1]$. The $L^2$ error is $7.01e-04$. All plots are obtained using physics-informed neural networks with hyper-parameters summarized in table \ref{tab: hyper-parameters}.}
        \label{fig: Stefan_1D2P_inverse_I}
\end{figure}

A visual assessment of the predictive accuracy of our framework is given in figure  \ref{fig: Stefan_1D2P_inverse_I}. In particular, the top panel shows a comparison against the exact temperature $u$ in the entire domain $\Omega$. A key observation is that the predicted solution achieves good accuracy with relative $L^2$ error $1.91e-03$, although there is a bit larger error along the two boundaries, which comes to no surprise because of the missing boundary conditions. Also, we present a more detailed assessment of the predicted interface in the bottom panel, which suggests that the proposed methodology is
able to correctly identify the unknown interface with good accuracy.

\subsubsection{Inverse one-dimensional two-phase Stefan problems: Type II}
\label{sec: Stefan_1D2P_inverse_II}

In this section, we employ the same example in section \ref{sec: Stefan_1D2P_direct} and \ref{sec: Stefan_1D2P_inverse_I}  but reformulate it as the second type of inverse Stefan problems, as described in section \ref{sec: Stefan_1D2P} and table \ref{tab: Stefan_1D2P}. To this end, we proceed by approximating $u_1(x,t)$, $u_2(x,t)$ and $s(t)$ with two independent fully connected neural networks which are denoted by $u_{\bm \theta}^{(1)}, u_{\bm \theta}^{(2)}$ and $s_{\bm \beta}$. The parameters of these networks can be calibrated by minimizing the sum of squared errors
\begin{align}
    \label{eq: loss_Stefan_1D2P_inverse_II}
    \mathcal{L}(\bm{\theta, \beta}) = \sum_{k=1}^2 \Big[\mathcal{L}_r^{(k)}(\bm{\theta}) + \mathcal{L}_{s_{bc}}^{(k)}(\bm{\theta, \beta})   \Big]  +  \mathcal{L}_{s_{Nc}}(\bm{\theta, \beta}) + \mathcal{L}_{s_0}(\bm{\beta}) + \mathcal{L}_{\text{data}}(\bm{\theta}),
\end{align}
where
\begin{align}
    \label{eq: loss_Stefan_1D2P_L_data}
    \mathcal{L}_{\text{data}}(\theta) = \frac{1}{M} \sum_{i=1}^{M} | u_{\bm{\theta}}(x^i_{\text{data}}, t^i_{\text{data}}) - u^i |^2.
\end{align}
The training data set $\{(x_{\text{data}}^i, t_{\text{data}}^i), u^i\}_{i=1}^M$ is generated by randomly sampling $M$ measurement points inside the domain $\Omega$ and obtain corresponding data for $u^i$ using equation \ref{eq: Stefan_1D2P_u_exact}. It is worth emphasizing that for any given data-point $\{(x_{\text{data}}^i, t_{\text{data}}^i), u^i)\}$, we do not know the corresponding equation that governs $u^i$ during training. 

\begin{figure}
     \centering
     \begin{subfigure}[b]{0.8\textwidth}
         \centering
         \includegraphics[width=\textwidth]{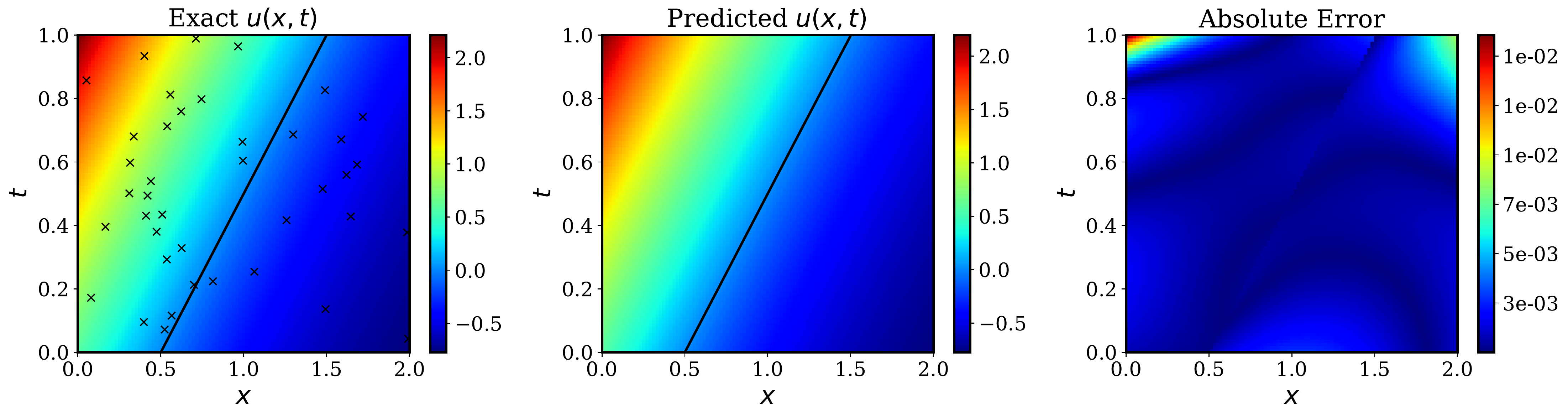}
         \caption{}
         \label{fig: Stefan_1D2P_inverse_II_pred_u}
     \end{subfigure}
     \begin{subfigure}[b]{0.5\textwidth}
         \centering
         \includegraphics[width=\textwidth]{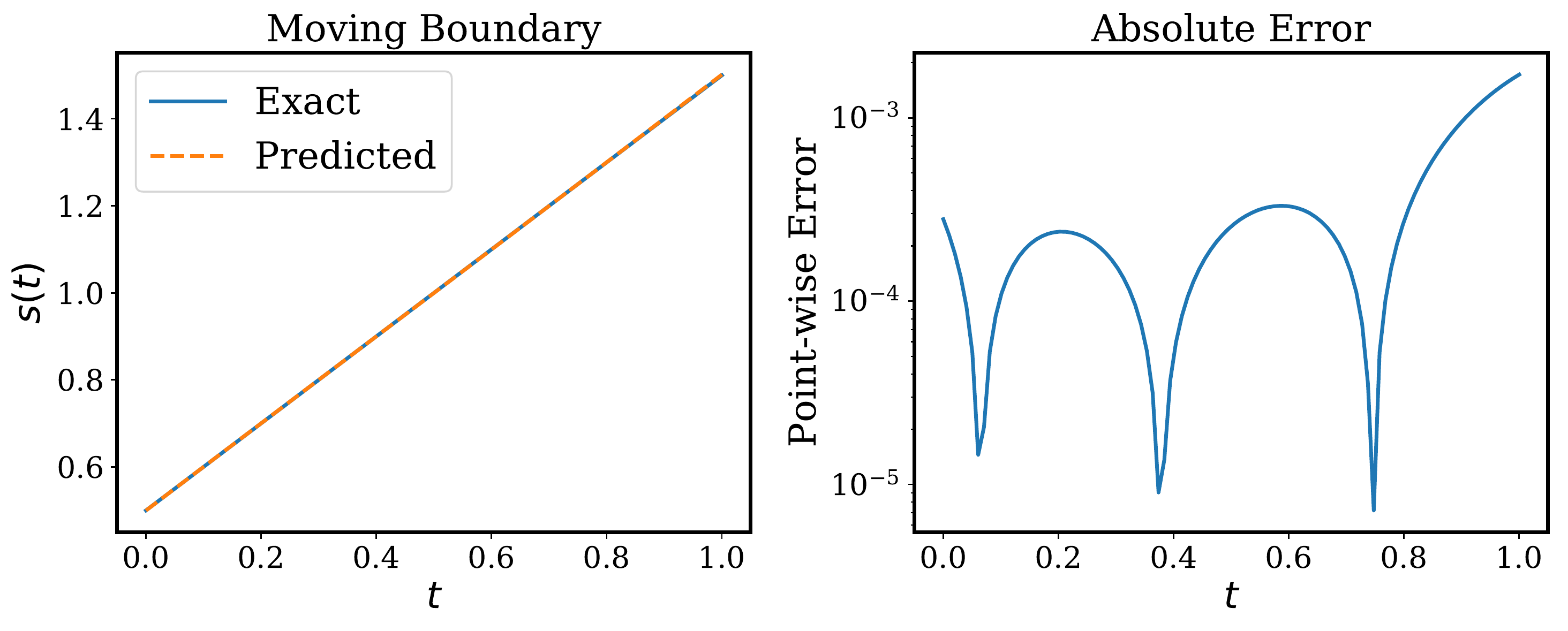}
         \caption{}
         \label{fig: Stefan_1D2P_inverse_II_pred_mb}
     \end{subfigure}
        \caption{{\em Inverse one-dimensional two-phase Stefan problem Type II: } (a) Exact temperature solution along with training data ($\times$) versus the predicted solution. The relative $L^2$ error: $2.57e-03$. (b) {\em Left:} Comparison of the and exact predicted moving boundary. {\em Right:} The absolute error between the exact  and the predicted free boundary for $t\in [0, 1]$. The $L^2$ error is $3.93e-04$. All plots are obtained using physics-informed neural networks with hyper-parameters summarized in table \ref{tab: hyper-parameters}.}
         \label{fig:Stefan_1D2P_inverse_II_pred}
\end{figure}

The data set alongside the predicted temperature solution as well as the interface are depicted in figure \ref{fig:Stefan_1D2P_inverse_II_pred}. These figures indicates that our algorithm is able to accurately learn the latent variables including $u_1(x,t), u_2(x,t)$ and $s(t)$ with a relative $L^2$ error of $2.57e-03$ and $3.93e-04$, respectively.

\begin{figure}
     \centering
     \begin{subfigure}[b]{0.3\textwidth}
         \centering
         \includegraphics[width=\textwidth]{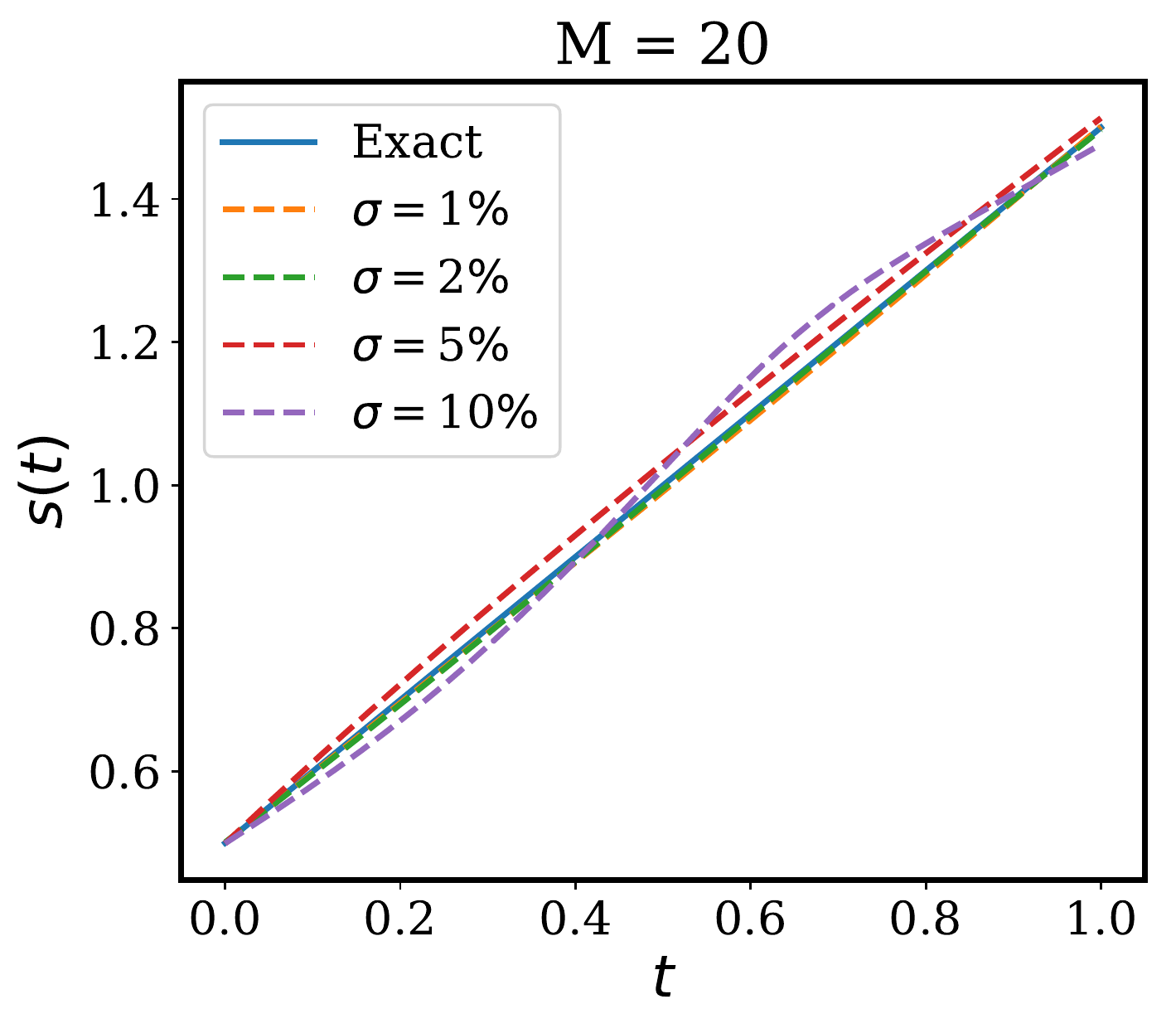}
         \caption{}
         \label{fig:Stefan_1D2P_inverse_II_pred_mb_M_20}
     \end{subfigure}
     \begin{subfigure}[b]{0.3\textwidth}
         \centering
         \includegraphics[width=\textwidth]{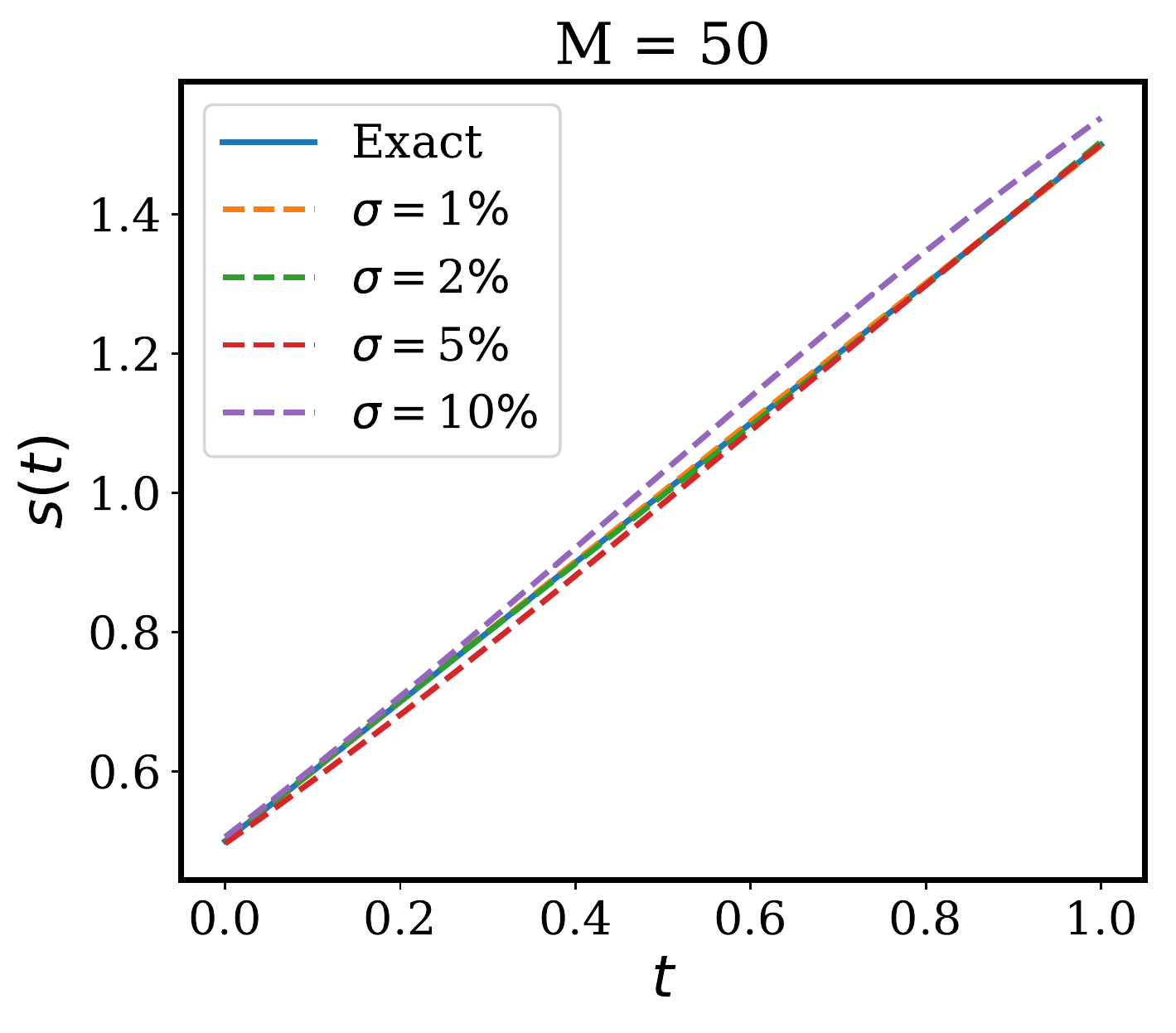}
         \caption{}
         \label{fig:Stefan_1D2P_inverse_II_pred_mb_M_50}
     \end{subfigure}
     \begin{subfigure}[b]{0.3\textwidth}
         \centering
         \includegraphics[width=\textwidth]{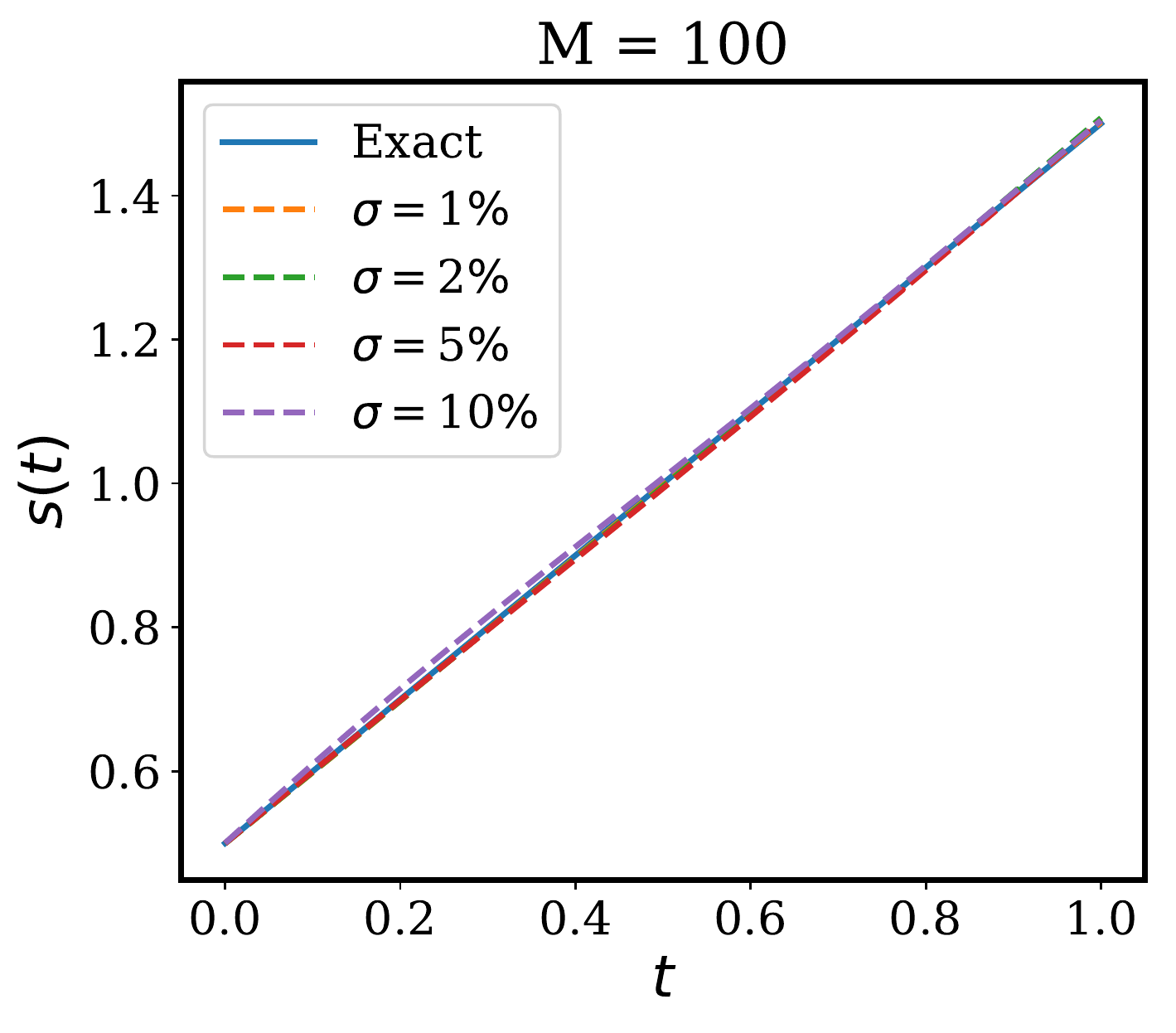}
         \caption{}
         \label{fig:Stefan_1D2P_inverse_II_pred_mb_M_100}
     \end{subfigure}
        \caption{{Inverse one-dimensional two-phase inverse Stefan problem Type II: } Comparison between the exact and the predicted moving boundary corresponding to different number of data points $M \in \{20, 50, 100\}$ and with different noise levels $\sigma \in \{1\%, 2\%, 5\%, 10\%\}$. All plots are obtained using physics-informed neural networks with hyper-parameters summarized in table \ref{tab: hyper-parameters}.}
        \label{fig:Stefan_1D2P_inverse_II_pred_mb_different_M}
\end{figure}

Furthermore, we performed a systematic study of the reported results in
figure \ref{fig:Stefan_1D2P_inverse_II_pred} with respect to noise levels and amount of data points by keeping the neural network architectures fixed to the settings described above. In particular, we added white noise with different magnitude according to the formula \ref{eq: noise} and varied the number of data points. The results of this study are summarized in tables \ref{tab:Stefan_2D1P_inverse_II_pred_u} and \ref{tab:Stefan_1D2P_inverse_II_pred_mb}. It can be concluded that using larger data set enhances the performance of our framework and mitigates the negative consequences of noise corruptions. Another fundamental point to make is that the predicted interface $s(t)$ obtained using our model remains accurate, even when the predicted solution is not very accurate.  This is also evident from a visual comparison of the predicted interface given in figure \ref{fig:Stefan_1D2P_inverse_II_pred_mb_different_M}.

\begin{table}[]
\renewcommand*{\arraystretch}{1.6}
    \centering
    \begin{tabular}{|c|c|c|c|c|c|}
\hline
\diagbox{Data Points M}{Noise $\sigma$}  & $\sigma =0\%$      & $\sigma =1\%$ &  $\sigma =2\%$ &  $\sigma =5\%$&  $\sigma =10\%$ \\ \hline
M = 20 &3.37e-03 & 2.89e-02    &  2.86e-02   &  9.22e-02  & 1.86e-01   \\ \hline
M = 50 &1.24e-03 & 6.58e-03     &   1.37e-02   &4.19e-02    & 8.40e-02   \\ \hline
M = 100&6.80e-03 & 6.37e-03    &1.37e-02   & 2.13e-02   & 5.17e-02     \\ \hline
\end{tabular}
    \caption{{\em Inverse one-dimensional two-phase Stefan problem Type II:} Relative $L^2$ errors of predicted solution with different number of data points and different noise level $\sigma$.}
    \label{tab:Stefan_1D2P_inverse_II_pred_u}
\end{table}

\begin{table}[]
\renewcommand*{\arraystretch}{1.6}
    \centering
    \begin{tabular}{|c|c|c|c|c|c|}
\hline
\diagbox{Data Points M}{Noise $\sigma$}  & $\sigma =0\%$      & $\sigma =1\%$ &  $\sigma =2\%$ &  $\sigma =5\%$&  $\sigma =10\%$ \\ \hline
M = 20 &8.07e-04 &  6.01e-03    &  4.37e-03 &  2.36e-02   &  3.04e-02    \\ \hline
M = 50 &3.75e-04 & 1.69e-03    & 2.24e-03  & 1.25e-02   & 3.05e-02     \\ \hline
M = 100&2.22e-04 & 1.96e-03   &  3.23e-03    & 4.26e-03   & 8.54e-03     \\ \hline
\end{tabular}
    \caption{{\em Inverse one-dimensional two-phase Stefan problem Type II:} Relative $L^2$ errors of predicted free boundary with different number of data points and different noise level $\sigma$.}
    \label{tab:Stefan_1D2P_inverse_II_pred_mb}
\end{table}

To test our algorithm even further, let us now take a closer look at heat transfer equation \ref{eq: Stefan_1D2P_eq} and ask: what would happen if the thermal diffusivity parameters $k_1, k_2$ are unknown? As we demonstrate in the following, the proposed framework is capable of handling such cases. Now our goal is to identify the unknown thermal diffusivities  $k_1$ and $k_2$, the unknown  interface $s(t)$, as well as to obtain an accurate reconstruction of the temperature distribution $u(x,t)$ in the whole domain $\Omega$. To this end, we still follow the same methodology we formulated in the previous subsections, albeit we now replace the original residual loss \ref{eq: loss_Stefan_1D2P_L_r} by

\begin{align}
    \mathcal{L}_r^{(k)}(\bm{\theta}) = \frac{1}{N} \sum_{i=1}^{N}|\frac{\partial u_{\bm \theta}^{(k)}}{\partial t}(x^i, t^i) -\lambda_k \frac{\partial^2 u_{\bm \theta}^{(k)}}{\partial x^2}(x^i, t^i)|^2, \quad k=1, 2,
\end{align}
where $\lambda_1, \lambda_2 $ are two extra trainable parameters. In addition, since we know that the thermal diffusivity cannot be negative, those two parameters are initialized by 0.1 and are constrained to remain positive during model training.
We hope that the shared parameters of the neural networks $u_{\bm \theta}, s_{\bm \beta}$ along with two parameters $\lambda_1, \lambda_2$ can be learned by minimizing this modified loss function. 

\begin{figure}
     \centering
     \begin{subfigure}[b]{0.4\textwidth}
         \centering
         \includegraphics[width=\textwidth]{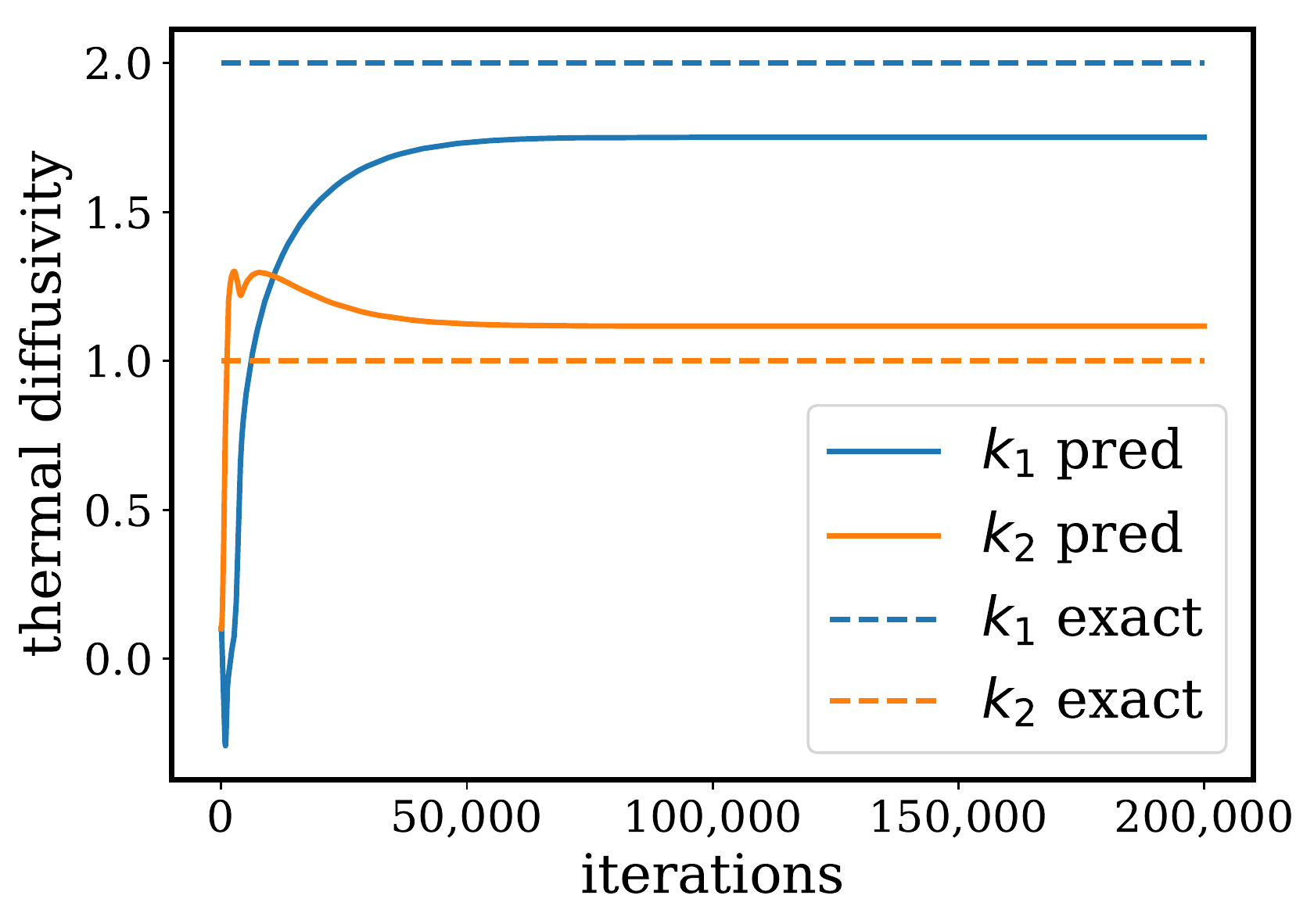}
         \caption{Inferred coefficients: $k_1=1.751$, $k_2=1.117$.}
         \label{fig: Stefan_1D2P_inverse_III_M1_k}
     \end{subfigure}
     \begin{subfigure}[b]{0.4\textwidth}
         \centering
         \includegraphics[width=\textwidth]{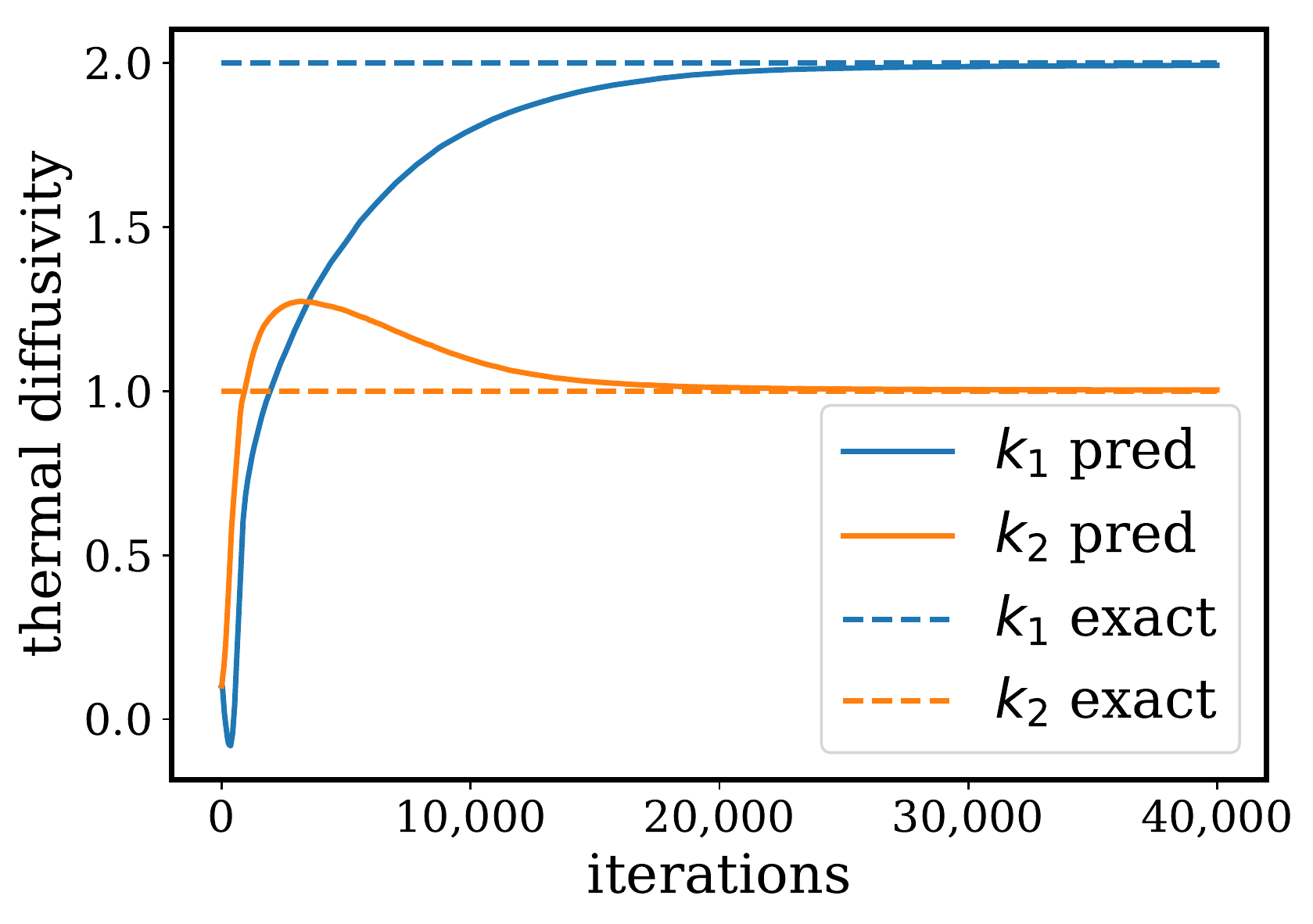}
         \caption{Inferred coefficients: $k_1 =1.993$, $k_2 =1.002$.}
         \label{fig: Stefan_1D2P_inverse_III_M2_k}
     \end{subfigure}
        \caption{{\em Inverse one-dimensional two-phase Stefan problem Type II with unknown parameters:} Convergence of the learned thermal diffusivity coefficients $k_1$ and $k_2$ using physics-informed neural networks with or without the adaptive learning rate annealing put forth by Wang {\em et. al.} \cite{wang2020understanding}.}
        \label{fig: Stefan_1D2P_inverse_III_M1_vs_M2}
\end{figure}

Since this task is considerably more complicated that the examples discussed in sections \ref{sec: Stefan_1D2P_direct} and \ref{sec: Stefan_1D2P_inverse_I}, we will employ a larger data set. In particular, the new data-set contains $M=200$ points of the temperature solution $u(x,t)$ which are randomly sampled in the whole domain $\Omega$. Figure \ref{fig: Stefan_1D2P_inverse_III_M1_k} shows the identified parameters $k_1$ and $k_2$ during the training of a standard physics-informed neural network \cite{raissi2019physics}. We can observe that the standard PINNs formulation of Raissi {\em et. al.} fails to correctly identify the unknown thermal diffusivity, even after $200,000$ training iterations. 
From the figure, it can be observed that two identified parameters do not change as training goes on. This observation suggests that our model seems to get stuck in some local minimum and as a result there is no hope to obtain better results just by increasing the training iterations under the current setting of hyper-parameters. Therefore, we point out that this result is not related to insufficient training iterations but to some issue pertaining to the model itself.  

To resolve this issue, let us apply the strategy of dynamics weights put forth in \cite{wang2020understanding} for automatically tuning the weights of the loss function \ref{eq: loss_Stefan_1D2P_inverse_II}. Specifically, let us  reformulate the loss function by
\begin{align}
    \label{eq: loss_Stefan_1D2P_inverse_II_adaptive}
    \mathcal{L}(\bm{\theta, \beta}) = \sum_{k=1}^2 \Big[\mathcal{L}_r^{(k)}(\bm{\theta}) + \mathcal{L}_{s_{bc}}^{(k)}(\bm{\theta, \beta})   \Big]  +  \mathcal{L}_{s_{Nc}}(\bm{\theta, \beta}) + \mathcal{L}_{s_0}(\bm{\beta}) + \lambda \mathcal{L}_{\text{data}}(\bm{\theta}),
\end{align}
where $\lambda$ is adaptively updated during training
by utilizing the back-propagated gradient statistics. 
Specifically, the estimates of $\lambda$ are computed by
\begin{align}
    \label{eq: esitmate_lambda}
    \hat{\lambda}^{(k+1)}=\max _{\bm \theta}\left\{\left|\nabla_{\theta} \mathcal{L}_{r}\right|\right\}  / \overline{\left|\nabla_{\theta} \lambda^{(k)} \mathcal{L}_{\text{data}}\right|},
\end{align}
where $\max _{\bm \theta}\left\{\left|\nabla_{\theta} \mathcal{L}_{r}\right|\right\}$ is the maximum value attained by $\left|\nabla_{\theta} \mathcal{L}_{r}\right|$ and  $\overline{\left|\nabla_{\theta} \lambda^{(k)} \mathcal{L}_{\text{data}}\right|}$ denotes the mean of $\left|\nabla_{\theta} \lambda^{(k)} \mathcal{L}_{\text{data}}\right|$  over the parameter space $\bm{\theta}$.
The weighting coefficients $\lambda$ for the next iteration are updated using a moving average of the form
\begin{align}
    \label{eq: moving_average}
    \lambda^{(k+1)}=(1-\alpha) \lambda^{(k)}+\alpha \hat{\lambda}^{(k+1)}
\end{align}
with $\alpha =0.1$. As shown in \cite{wang2020understanding}, this adaptive strategy can effectively mitigate pathologies arising in the training of physics-informed neural networks due to stiffness in gradient flow dynamics.

Figure \ref{fig: Stefan_1D2P_inverse_III_M1_vs_M2} and table  \ref{tab: Stefan_1D2P_inverse_III_parameters} present comparisons of the identified parameters between the original PINNs formulation of Raissi {\em et. al.} \cite{raissi2019physics}, and the proposed PINNs formulation with adaptive weights \cite{wang2020understanding}. Notice that the PINNs with adaptive weights not only converges to the exact parameters much faster, but also yields a significantly improved identification accuracy. In addition, we also investigate the accuracy of our  reconstructed temperature solution $u(x,t)$ and inferred interface $s(t)$ with respect to these two methods.
A comparison of relative $L^2$ error in $u(x,t)$ and $s(t)$ between these two models is presented in table \ref{tab: Stefan_1D2P_inverse_III_error} from which we can see that the dynamic weights approach improves the relative prediction error by about one order of magnitude.  These figures and tables highly suggest that the weights in the loss function play a quite important role, and  choosing appropriate weighting coefficients can enhance the performance of PINNs by accelerating convergence and avoiding bad local minima.

\begin{table}
\renewcommand*{\arraystretch}{1.8}
    \centering
   \begin{tabular}{|c|c|c|}
\hline
Correct PDE               & $\frac{\partial u_1}{\partial t} - 2\frac{\partial^2 u_1}{\partial x^2} = 0$ & $\frac{\partial u_1}{\partial t} - \frac{\partial^2 u_1}{\partial x^2} = 0$ \\ \hline
Identified PDE (Original) & $\frac{\partial u_1}{\partial t} - 1.712\frac{\partial^2 u_1}{\partial x^2} = 0$  &$\frac{\partial u_1}{\partial t} - 1.137\frac{\partial^2 u_1}{\partial x^2} = 0$ \\ \hline
Identified PDE (Adaptive)  &$\frac{\partial u_1}{\partial t} - 1.993\frac{\partial^2 u_1}{\partial x^2} = 0$      &  $\frac{\partial u_1}{\partial t} - 1.002\frac{\partial^2 u_1}{\partial x^2} = 0$ \\ \hline
\end{tabular}
    \caption{{\em Inverse one-dimensional two-phase Stefan problem Type II with unknown parameters:} Correct partial differential equation along with the identified one obtained using physics-informed neural networks with or without the adaptive learning rate annealing put forth by Wang {\em et. al.} \cite{wang2020understanding}.}
    \label{tab: Stefan_1D2P_inverse_III_parameters}
\end{table}

\begin{table}
    \centering
    \renewcommand*{\arraystretch}{1.6}
\begin{tabular}{|c|c|c|}
\hline
\diagbox{Method}{$L^2$ error}     & $L^2$ error in $u$ & $L^2$ error in $s$\\ \hline
Original & 1.01e-02 & 4.28e-03 \\ \hline
Adaptive & 1.03e-03 & 4.36e-04 \\ \hline
\end{tabular}
    \caption{{\em Inverse one-dimensional two-phase Stefan problem Type II with unknown parameters:} Relative $L^2$ errors for the predicted solution and free boundary using physics-informed neural networks with or without the adaptive learning rate annealing put forth by Wang {\em et. al.} \cite{wang2020understanding}.}
    \label{tab: Stefan_1D2P_inverse_III_error}
\end{table}

\subsection{Two-dimensional one-phase Stefan problems}
\label{sec: Stefan_2D1P}

As our last example, we want to examine the validity of our model for two-dimensional Stefan problems. Compared to the case of one space variable in which there exist numerous results concerning the existence, uniqueness and numerical computation of the solution, in the case of two or more space variables, the literature is rather scarce, see \cite{nochetto1991adaptive, bonnerot1977numerical,beckett2001moving, colton1974inverse,jochum1982numerical,colton1984numerical}. One main reason is that the geometry is particularly simple in one-dimensional case, while can be extremely complicated in  high-dimensional cases, not to mention the formidable mathematical and numerical difficulties present.

Let us start by giving a general formulation for the high-dimensional single phase Stefan problem. For $0 < t \leq T < \infty$, let $\Omega = \bigcup_{0 < t \leq T} \{(\bm{x},t) \in \R^{n+1} : x \in \Omega(t)\}$ be a bounded connected domain in $\R^n$ with boundary $\partial \Omega(t) = \Gamma(t) \bigcup \Sigma(t)$, where $\Gamma = \cup_{0 < t <T} \Gamma(t)$ denotes the ``fixed" boundaries and $\Sigma = \cup_{0 < t <T} \Sigma(t)$ denotes the ``free" boundaries. Here we assume that $\Gamma$ is simply connected whereas $\Sigma$ may be multiply connected. In particular, $\bm{x} = (x_1, \cdots, x_n)$ will denote a point in Euclidean space $\R^n$. The equations describing the Stefan problem are as follows:
\begin{align}
    \label{eq: Stefan_2D1P_eq}
    &u_t - \Delta u =0 \quad \text{in } \Omega \\
     \label{eq: Stefan_2D1P_bc}
    &u|_{\overline{\Gamma}} = g \\
    \label{eq: Stefan_2D1P_ic}
    &u|_{\Omega(0)} = u_0  \\
      \label{eq: Stefan_2D1P_Sbc}
    &u|_{\overline{\Sigma}} = u^* \\
    \label{eq: Stefan_2D1P_Nc}
    & \frac{\partial u}{\partial n}\Big|_{\overline{\Sigma}} = h \\
        \label{eq: Stefan_2D1P_s0}
    &    \Sigma(0) = \Sigma_0,
\end{align}
where $\Delta$ is the Laplace operator with respect to spatial variables and $n$ is the normal vector with respect to the space variables in the domain $\Omega$. For direct Stefan problem, the goal would be to determine the unknown free surface $\Sigma(t)$  as well as the temperature solution $u(\bm{x},t)$ satisfying equations \ref{eq: Stefan_2D1P_eq} - \ref{eq: Stefan_2D1P_ic}. 

For inverse two-dimensional Stefan problems, instead of assuming knowledge of the boundary conditions, additional information at the final time $T$ is provided
\begin{align}
\label{eq: Stefan_2D1P_final_time}
    u|_{\Omega(T)} = u_T
\end{align}
The objective of this class of inverse Stefan problem, which we refer to as ``Inverse type I", is to find the temperature solution $u(\bm{x},t)$ and the free surface $\Sigma(t)$ satisfying equations \ref{eq: Stefan_1D2P_eq}, \ref{eq: Stefan_2D1P_ic}, \ref{eq: Stefan_2D1P_Sbc} - \ref{eq: Stefan_2D1P_final_time}. 

As in the first two case studies, the second type inverse Stefan problem which we refer to as ``Inverse type II" is to recover the temperature distribution $u(\bm{x},t)$ as well as identify the free boundary $\Sigma(t)$, assuming that some measurements of the temperature inside the domain are provided. A summary of these three Stefan type problems with known conditions and objective is presented in table \ref{tab: Stefan_2D1P}. 

 \begin{table}
    \centering
     \renewcommand*{\arraystretch}{1.8}
   \begin{tabular}{|c|c|c|c|}
\hline
                & Equations  & Observed & Latent \\ \hline
Direct          & \ref{eq: Stefan_2D1P_eq} - \ref{eq: Stefan_2D1P_ic}       &  $u_0, g,h, \Sigma_0 $ & $u(\bm{x},t), \Sigma(t)$       \\ \hline
Inverse Type I  & \ref{eq: Stefan_1D2P_eq}, \ref{eq: Stefan_2D1P_ic}, \ref{eq: Stefan_2D1P_Sbc} - \ref{eq: Stefan_2D1P_final_time}      & $c, u_0, u_T,g, h, \Sigma_0$ & $u(\bm{x},t), \Sigma(t) $     \\ \hline
Inverse Type II &\ref{eq: Stefan_2D1P_eq},  \ref{eq: Stefan_2D1P_Sbc} - \ref{eq: Stefan_1D12P_Sic},       & $\{(\bm{x}^j_{\text{data}}, t^j_{\text{data}}), u^j \}_{j=1}^M, g, h, \Sigma_0 $    &  $u(\bm{x},t), \Sigma(t)$   \\ \hline
\end{tabular}
    \caption{{\em Two-dimensional one-phase Stefan problems:} Summary of conditions and the objective of each type of Stefan problem formulated in section \ref{sec: Stefan_2D1P}.}
    \label{tab: Stefan_2D1P}
\end{table}

\subsubsection{Direct two-dimensional one-phase Stefan problems}
\label{sec: Stefan_2D1P_direct}
We first consider a specific direct two-dimensional single phase Stefan example \cite{colton1984numerical}. Let $0 \leq t \leq T =1$ and the computational domain given by 
\begin{align}
    \Omega(t) = \{(x, y, t) \in \R^3 : 0 < x < s(y,t), 0< y <1\} \subset \Omega^*,
\end{align}
where $s(y ,t)$ denotes the unknown free surface and $\Omega^* = [0,2.25] \times [0,1] \times [0,1]$. Let
\begin{align}
    \Phi(x, y, t) = x - s(y,t),
\end{align}
and define $\Sigma(t)$ by

\begin{align}
    \Sigma(t) = \{(x, y, t) \in \R^3 : \Phi(x, y, t) = 0, 0< y < 1, 0 \leq t \leq 1\}.
\end{align}
As described in table \ref{tab: Stefan_2D1P}, we aim to solve equations \ref{eq: Stefan_2D1P_eq} - \ref{eq: Stefan_2D1P_s0} with 
\begin{align}
    &u_0(x, y) = \exp(-x + \frac{1}{2}y + \frac{1}{2}),  \quad x,y \in \Omega(0) \\
    &u(s(y,t), y, t) = u^* = 0 \\
    &  h =  \frac{1}{|\nabla \Phi|}\frac{\partial \Phi}{\partial t} \\
   & s(y,0) = s_0(y) = \frac{1}{2}y + \frac{1}{2} .
\end{align}
The boundary conditions considered in the direct problem are given by
\begin{align}
    \label{eq: Stefan_2D_bc1}
    &u(x,0,t) = g_1(x, t) = \exp(1.25 t - x + 1/2) -1,  \quad (x,0,t) \in \Omega \\
     \label{eq: Stefan_2D_bc2}
    &u(0, y, t) = g_2(y,t) = 1.25  t + 0.5 y + 1/2) - 1, \quad (0, y, t) \in \Omega \\
    &u(x,1,t) = g_3(x,t) = \exp(1.25t - x + 1), \quad (x,1,t) \in \Omega.
\end{align}
The exact solution of this problem is given by
\begin{align}
    \label{eq: Stefan_2D_u_sol}
    &u(x,y,t) = \exp(\frac{5}{4}t - x + \frac{1}{2}y + \frac{1}{2}) -1 \\
    \label{eq: Stefan_2D_s_sol}
    &s(y, t) = \frac{1}{2}y + \frac{5}{4}t + \frac{1}{2}.
\end{align}

Similarly, we represent the temperature distribution $u(x,y,t)$ and the free surface $s(y,t)$ by two separate deep fully-connected neural networks. These  networks are trained by minimizing the sum of squared errors loss of equation 
 \begin{align}
     \label{eq: loss_Stefan_2D1P_direct}
     \mathcal{L}(\bm{\theta, \beta}) =   \mathcal{L}_r(\bm{\theta}) +  \mathcal{L}_{u_0}(\bm{\theta}) +  \mathcal{L}_{u_{bc}}(\bm{\theta})  +  \mathcal{L}_{S_{bc}}(\bm{\theta, \beta}) +  \mathcal{L}_{s_{Nc}}(\bm{\theta, \beta}) + \mathcal{L}_{s_0}(\bm{\beta}), 
 \end{align}
where 
\begin{align}
    \label{eq: loss_Stefan_2D1P_L_r}
    &\mathcal{L}_r(\bm{\theta}) = \frac{1}{N} \sum_{i=1}^N |u_{\bm{\theta}}(x^i, y^i, t^i) - u(x^i, y^i, t^i)|^2 \\
    \label{eq: loss_Stefan_2D1P_L_u0}
    & \mathcal{L}_{u_0}(\theta) = \frac{1}{N} \sum_{i=1}^N |u_{\bm{\theta}}(x^i, y^i, 0) - u(x^i, y^i, 0) |^2 \\
    \label{eq: loss_Stefan_2D1P_L_bc}
    & \mathcal{L}_{u_{bc}}(\bm{\theta}) = \frac{1}{N} \sum_{i=1}^N \Big[ |u_{\bm{\theta}}(x^i, 0, t^i) - g_1(x^i, t^i) |^2 + |u_\theta(0, y^i, t^i) - g_2(y^i, t^i) |^2 +  |u_\theta(x^i, 1, t^i) - g_2(x^i, t^i) |^2 \Big]  \\
    \label{eq: loss_Stefan_2D1P_L_sbc}
    & \mathcal{L}_{s_{bc}}(\bm{\theta, \beta}) = \frac{1}{N} \sum_{i=1}^{N} |u_{\bm{\theta}}(s_{\bm{\beta}}(y^i, t^i), y^i, t^i)|^2 \\
     \label{eq: loss_Stefan_2D_L_sNc}
    &  \mathcal{L}_{s_{Nc}}(\bm{\theta, \beta}) = \frac{1}{N} \sum_{i=1}^{N} |\frac{\partial u_{\bm{\theta}}}{\partial \bm{n}}(s_{\bm{ \beta}}(y^i,t^i),y^i t^i) - h(y^i, t^i)|^2 \\
     \label{eq: loss_Stefan_2D_L_s0}
    & \mathcal{L}_{s_0}(\bm{ \beta})  =  \frac{1}{N} \sum_{i=1}^{N} |s_{\bm{\beta}}(y^i, 0) - s_0(y^i)|^2.
\end{align}
where $N$ is the batch size. It should be emphasised again that each mini-batch $\{(x^i,y^i, t^i)\}_{i=1}^N$ is randomly sampled from the artificial domain $\Omega^*$ rather than $\Omega$ because the free surface $s(y,t)$ is unknown. Then the predicted solution $u(x,y,t)$ is obtained by restricting to the domain $\Omega$ using the predicted free surface $s(y,t)$. 

\begin{figure}
     \centering
     \begin{subfigure}[b]{0.8\textwidth}
         \centering
         \includegraphics[width=\textwidth]{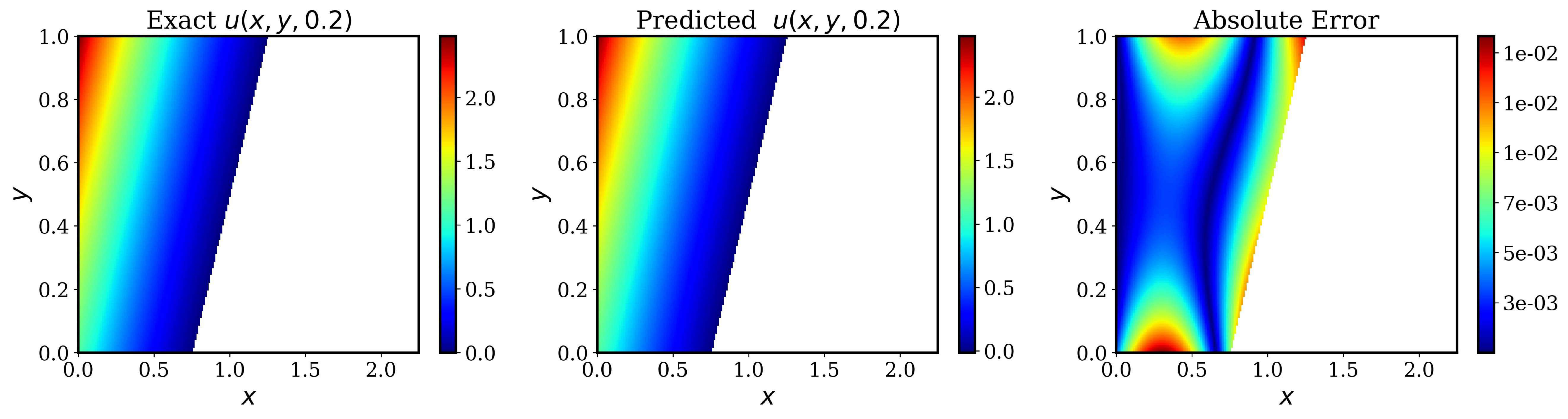}
     \end{subfigure}
     \begin{subfigure}[b]{0.8\textwidth}
         \centering
         \includegraphics[width=\textwidth]{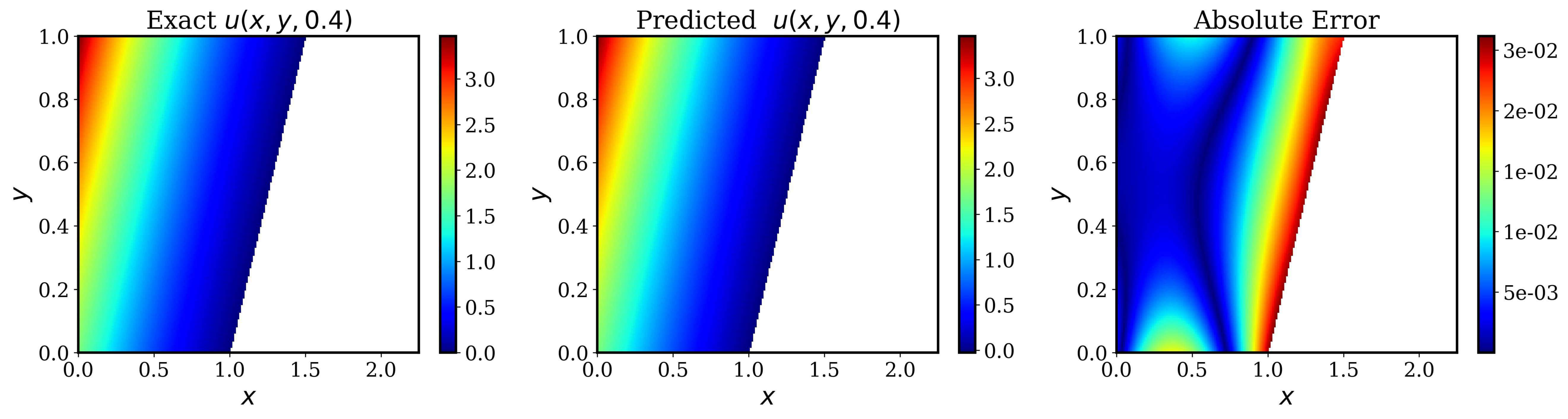}
     \end{subfigure}
     \begin{subfigure}[b]{0.8\textwidth}
         \centering
         \includegraphics[width=\textwidth]{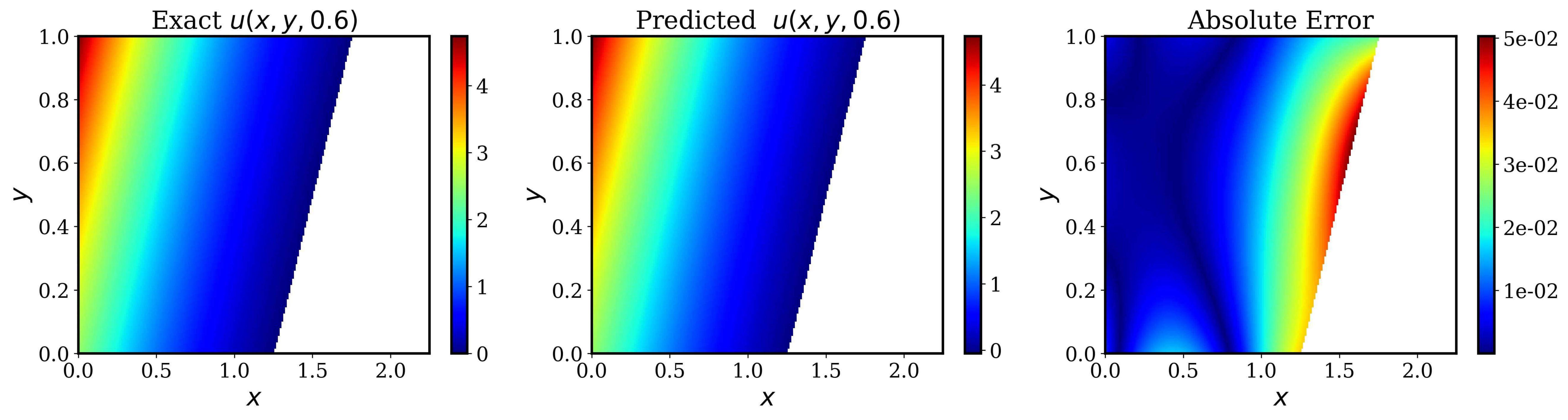}
     \end{subfigure}
      \begin{subfigure}[b]{0.8\textwidth}
         \centering
         \includegraphics[width=\textwidth]{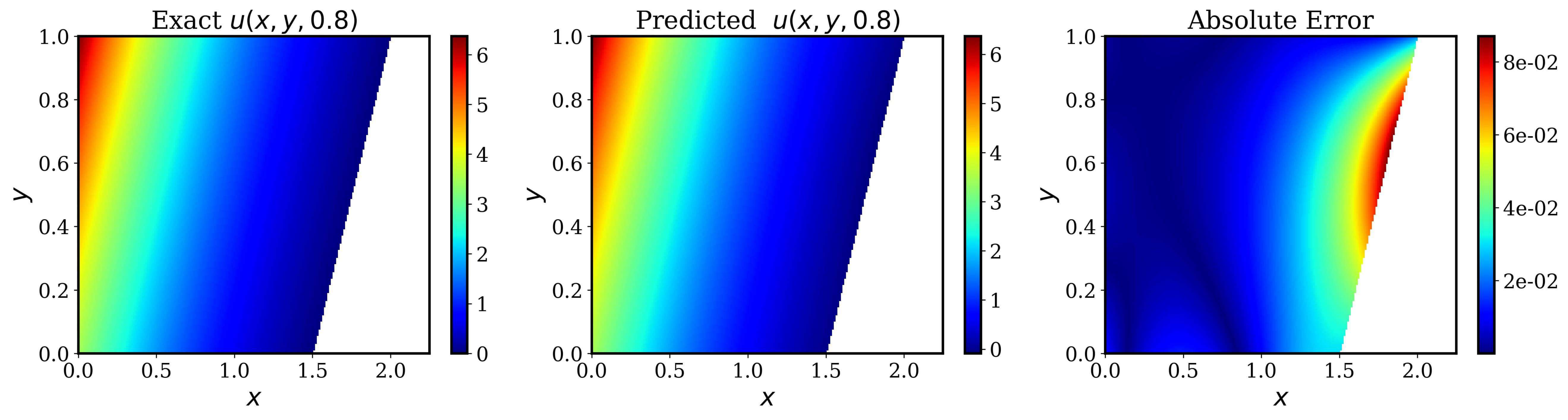}
     \end{subfigure}
        \caption{{\em Direct two-dimensional one-phase Stefan problem :} Comparison of the predicted and exact solutions corresponding to four different temporal snapshots.}
        \label{fig: Stefan2D_direct_snapshots}
\end{figure}

\begin{figure}
    \centering
    \includegraphics[width=0.8\textwidth]{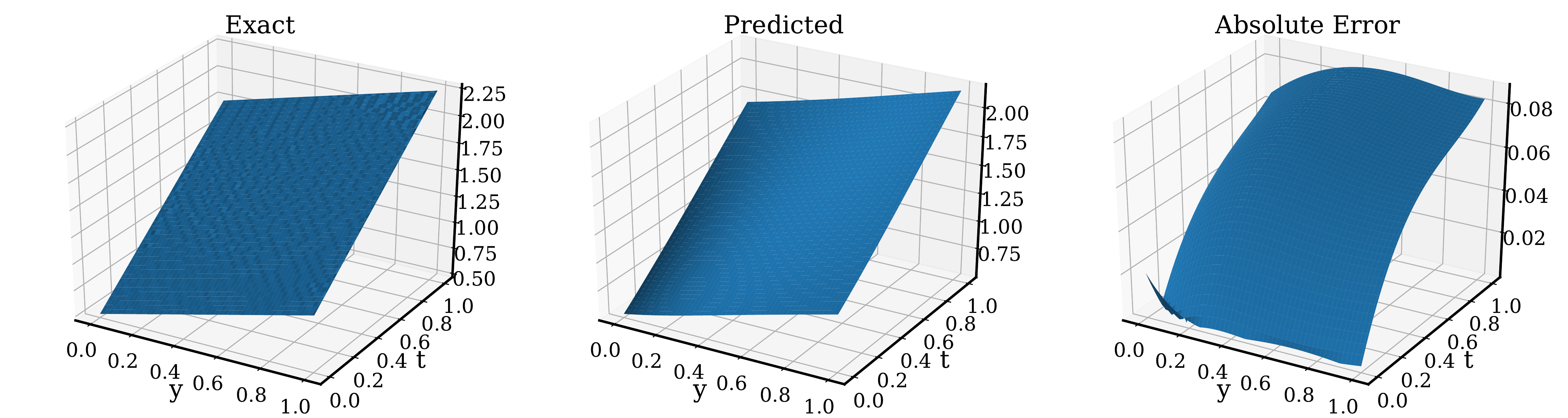}
    \caption{{\em Direct two-dimensional one-phase Stefan problem :}  Comparison of the exact and predicted free surface. The relative $L^2$ error is 4.32e-02.}
    \label{fig: Stefan_2D1P_direct_mb}
\end{figure}

Four representative snapshots of the exact and the predicted temperature distributions $u(x,y,t)$ are shown in figure \ref{fig: Stefan2D_direct_snapshots}. These snapshots are $\Delta t = 0.2$ apart and stretch from time $t = 0.2$ to $t = 0.8$. As it can be seen, a good agreement can be achieved between the predictions of our model and the exact solutions and the $L^\infty$ error we obtain is of order $O(10^{-2})$. In particular,  we compare the predicted free surface $s(y,t)$ to the exact one with the relative $L^2$ error $4.32-02$, as presented in figure \ref{fig: Stefan_2D1P_direct_mb}. 

\subsubsection{Inverse two-dimensional one-phase Stefan problems: Type I}

In section, we reformulate the example in \ref{sec: Stefan_2D1P_direct} as the first type of inverse two-dimensional Stefan problem whose conditions and objectives are summarized in table \ref{tab: Stefan_2D1P}. Recall that we seek to find the solution $u(x, y, t)$ and free boundary $s(y,t)$ satisfying \ref{eq: Stefan_2D1P_eq}, \ref{eq: Stefan_2D1P_ic}, \ref{eq: Stefan_2D1P_Sbc} - \ref{eq: Stefan_2D1P_s0} supplemented with the data at the final time $T=1$,
\begin{align}
    u(x, y, T) = u_T = \exp(1.25  - x + 0.5 y + 0.5) - 1.
\end{align}
To this end, we approximate the temperature distribution $u(x,y,t)$ and the free surface $s(y,t)$ by two separate neural networks $u_{\bm \theta}$ and $s_{\bm \beta}$. The shared parameters of the neural networks can be learned by  minimizing the following sum of squared errors loss function
\begin{align}
    \label{eq: loss_Stefan_2D_inverse_I}
     \mathcal{L}(\bm{\theta, \beta}) =   \mathcal{L}_r(\bm{\theta}) +  \mathcal{L}_{u_0}(\bm{\theta}) +  \mathcal{L}_{u_{T}}(\bm{\theta})  +  \mathcal{L}_{s_{bc}}(\bm{\theta, \beta}) +  \mathcal{L}_{s_{Nc}}(\bm{\theta, \beta}) + \mathcal{L}_{s_0}(\bm{\beta}), 
\end{align}
where all definitions of summands are exactly same as in the direct 2D Stefan problem except $\mathcal{L}_{u_{T}}(\theta) $, which is defined by the mean square loss function
\begin{align}
    \label{eq: loss_Stefan_2D_uT}
    \mathcal{L}_{u_T}(\bm{\theta}) = \frac{1}{N} \sum_{i=1}^N |u_{\bm{\theta}}(x^i, y^i, 1) - u(x^i, y^i, 1) |^2. \\
\end{align}

\begin{figure}
     \centering
      \begin{subfigure}[b]{0.8\textwidth}
         \centering
         \includegraphics[width=\textwidth]{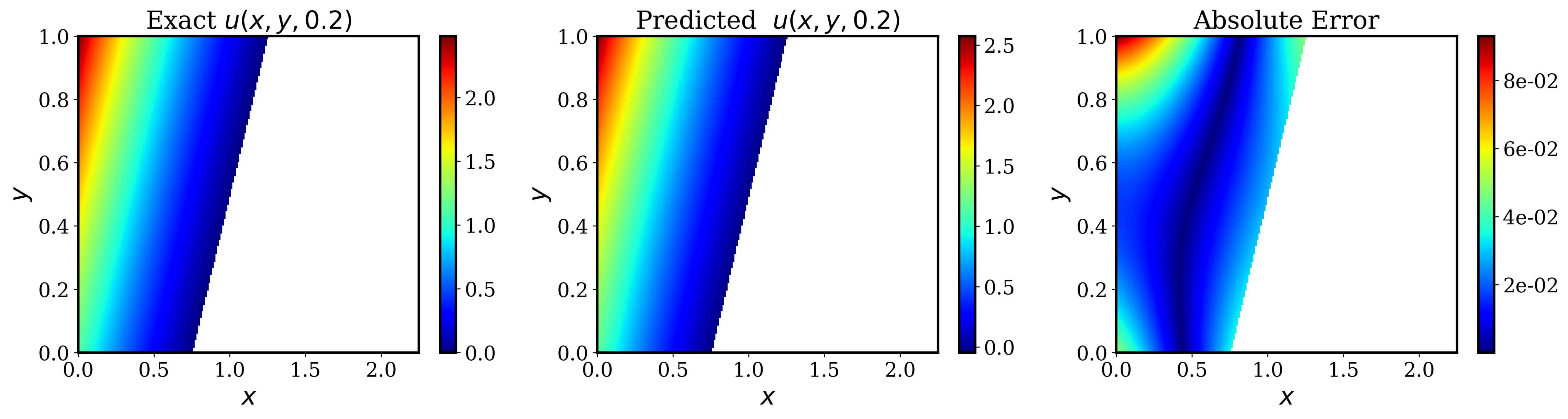}
     \end{subfigure}
     \begin{subfigure}[b]{0.8\textwidth}
         \centering
         \includegraphics[width=\textwidth]{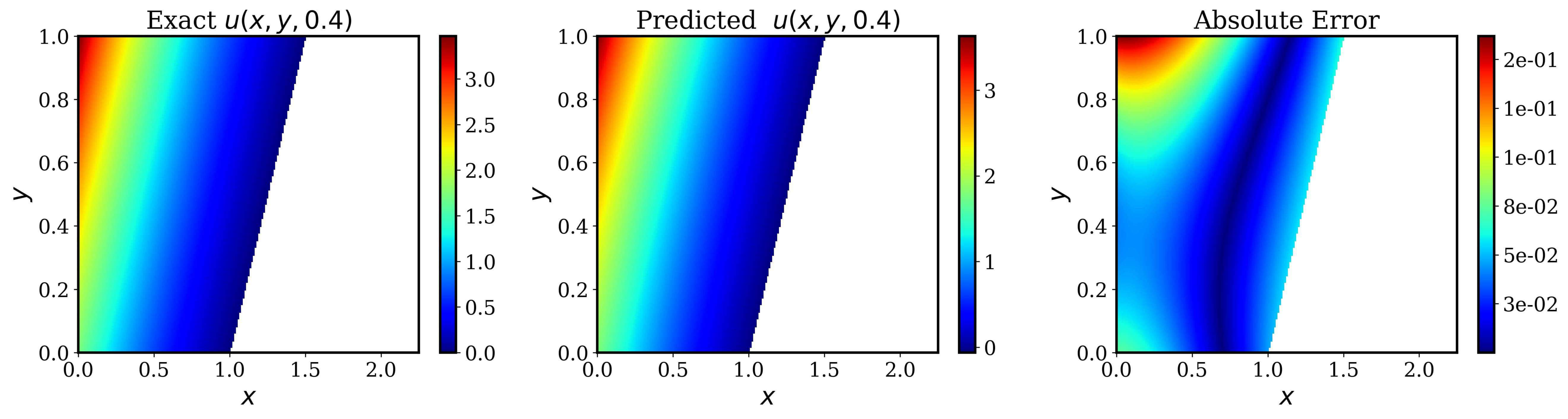}
     \end{subfigure}
     \begin{subfigure}[b]{0.8\textwidth}
         \centering
         \includegraphics[width=\textwidth]{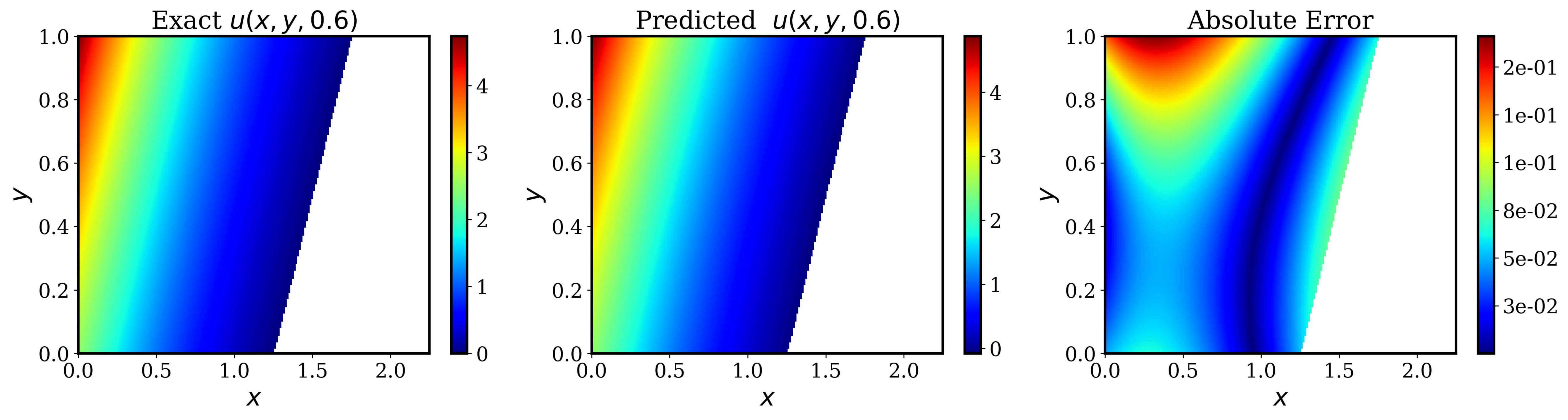}
     \end{subfigure}
      \begin{subfigure}[b]{0.8\textwidth}
         \centering
         \includegraphics[width=\textwidth]{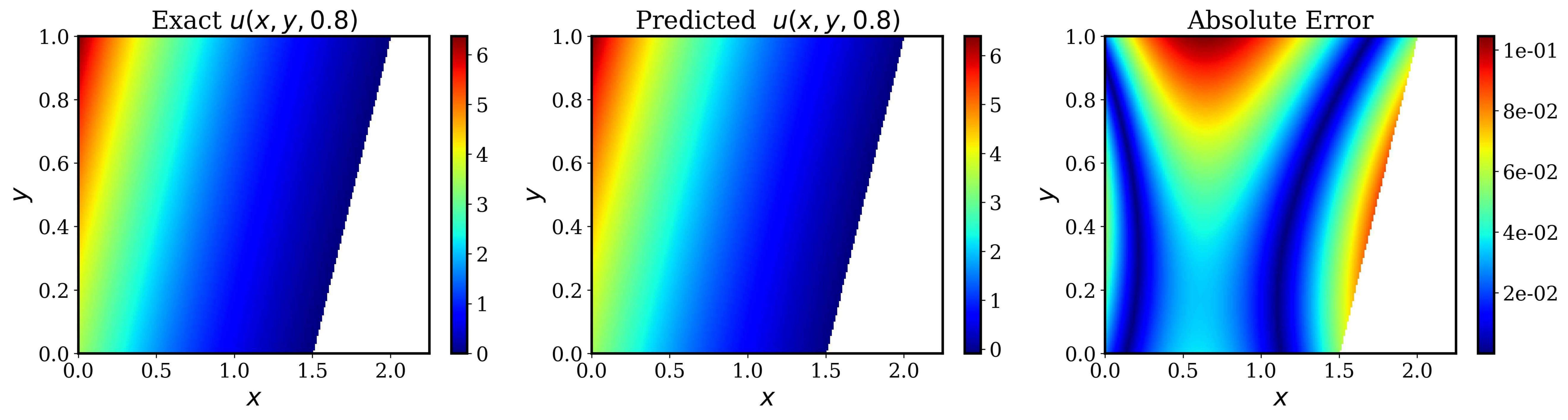}
     \end{subfigure}
        \caption{{\em Inverse two-dimensional one-phase Stefan problem Type I:} Comparison of the predicted and exact
         solutions corresponding to the four temporal snapshots.}
        \label{fig: Stefan2D_inverse_I_snapshots}
\end{figure}

Figure \ref{fig: Stefan2D_inverse_I_snapshots} depicts four snapshots of predicted temperature solution $u(x,y,t)$ as
well as the exact solution at different time instants
$t = 0.2, 0.4, 0.6, 0.8$, observing a good agreement between the two. In addition, a detailed comparison of the identified and the exact free boundary $s(y,t)$ is shown in figure \ref{fig: Stefan_1D2P_inverse_I_pred_mb} with a 
 resulting relative error of $4\%$. Some discrepancy, however, can be observed along the $t$ axis, which can be attributed to the lack of data for the free surface at the final time.

\begin{figure}
    \centering
    \includegraphics[width=0.8\textwidth]{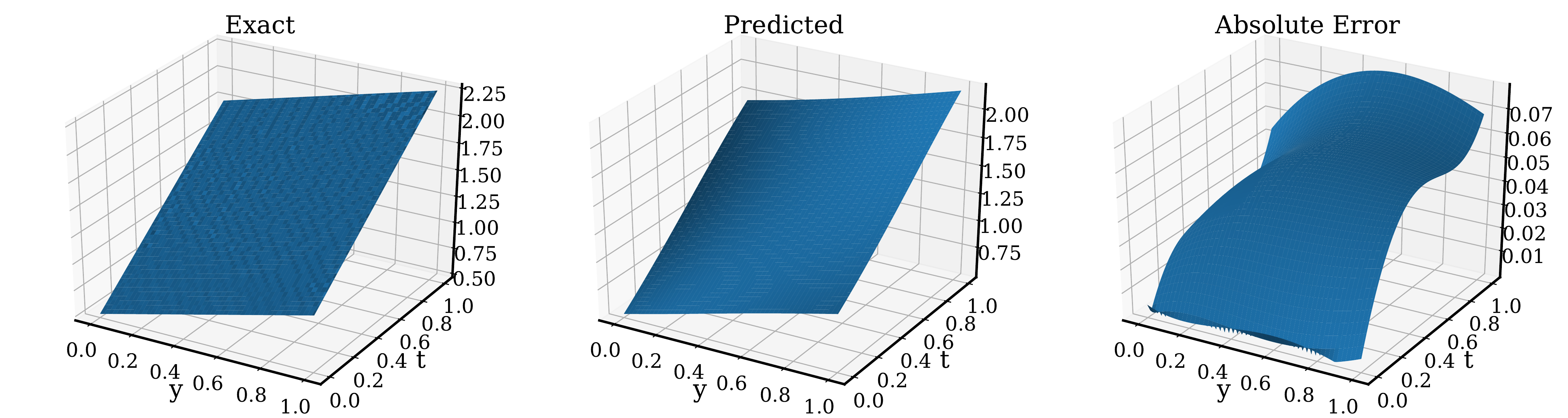}
    \caption{{\em Inverse two-dimensional one-phase Stefan problem Type I :} Comparison of the exact and predicted free surface. The relative $L^2$ error is 3.41e-02.}
    \label{fig:Stefan2D_inverse_I_mb}
\end{figure}

\subsubsection{Inverse two-dimensional one-phase Stefan problems: Type II}

As our last example, we consider the second type of inverse Stefan problem described in section \ref{sec: Stefan_2D1P} and table \ref{tab: Stefan_2D1P}. Here we emphasize that the free boundary $s(y,t)$ is unknown and no initial or boundary data is given. As before, we use two different neural networks to  represent the solution $u(x,y,t)$ and $s(y,t)$, denoted by $u_{\bm \theta}$ and $s_{\bm \beta}$ respectively. And then we can obtain the learned temperature distribution along with the free boundary by minimizing the following loss function
\begin{align}
    \label{eq: loss_Stefan_2D_inverse_II}
      \mathcal{L}(\bm{\theta, \beta}) =  \mathcal{L}_{\text{data}}(\bm{\theta}) +  \mathcal{L}_r(\bm{\theta}) +  \mathcal{L}_{s_{bc}}(\bm{\theta, \beta}) +  \mathcal{L}_{s_{Nc}}(\bm{\theta, \beta})  +   \mathcal{L}_{s_0}(\bm{\beta}),
\end{align}
where $\mathcal{L}_r(\bm{\theta}), \mathcal{L}_{s_{bc}}(\bm{\theta, \beta}), \mathcal{L}_{s_{Nc}}(\bm{\theta, \beta}) ,  \mathcal{L}_{s_0}(\bm{\beta})$ are exactly the same as in equations \ref{eq: loss_Stefan_2D1P_L_r}, \ref{eq: loss_Stefan_2D1P_L_sbc} - \ref{eq: loss_Stefan_2D_L_s0}, while  $\mathcal{L}_{\text{data}}$ is defined by
\begin{align}
    \label{eq: loss_Stefan_2D1P_L_data}
    \mathcal{L}_{\text{data}}(\bm{\theta})  = \frac{1}{M} \sum_{i=1}^M |u_{\bm \theta}(x^i_{\text{data}}, y^i_{\text{data}}, t^i_{\text{data}}) - u^i  |^2,
\end{align}
where $N$ is the batch size and $\{(x^i_{\text{data}}, y^i_{\text{data}}, t^i_{\text{data}}), u^i\}_{i=1}^M$ is the training data which is randomly sampled in the domain $\Omega$. 

\begin{figure}
     \centering
      \begin{subfigure}[b]{0.8\textwidth}
         \centering
         \includegraphics[width=\textwidth]{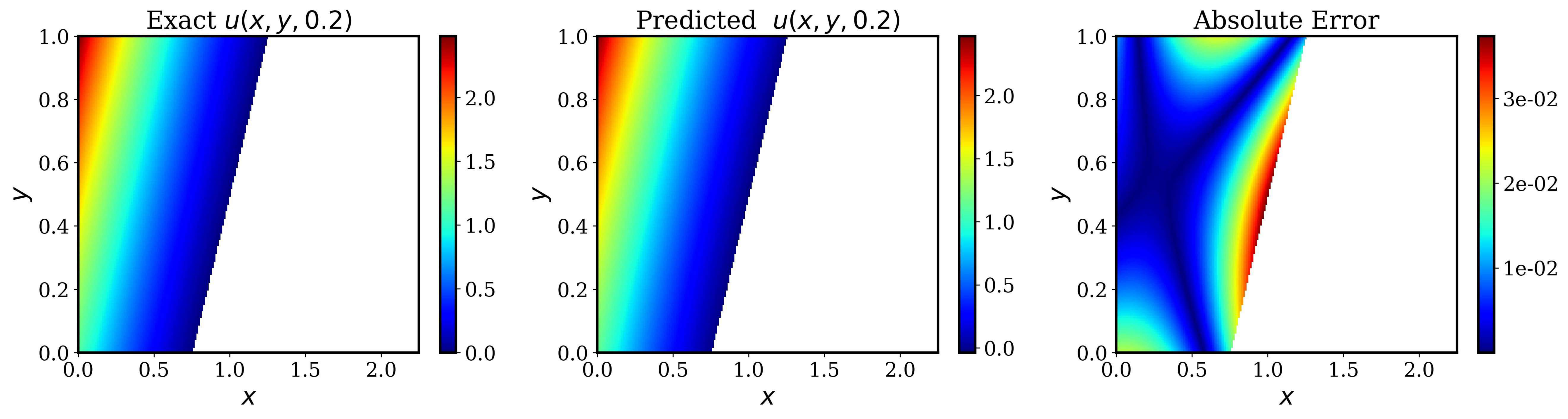}
     \end{subfigure}
     \begin{subfigure}[b]{0.8\textwidth}
         \centering
         \includegraphics[width=\textwidth]{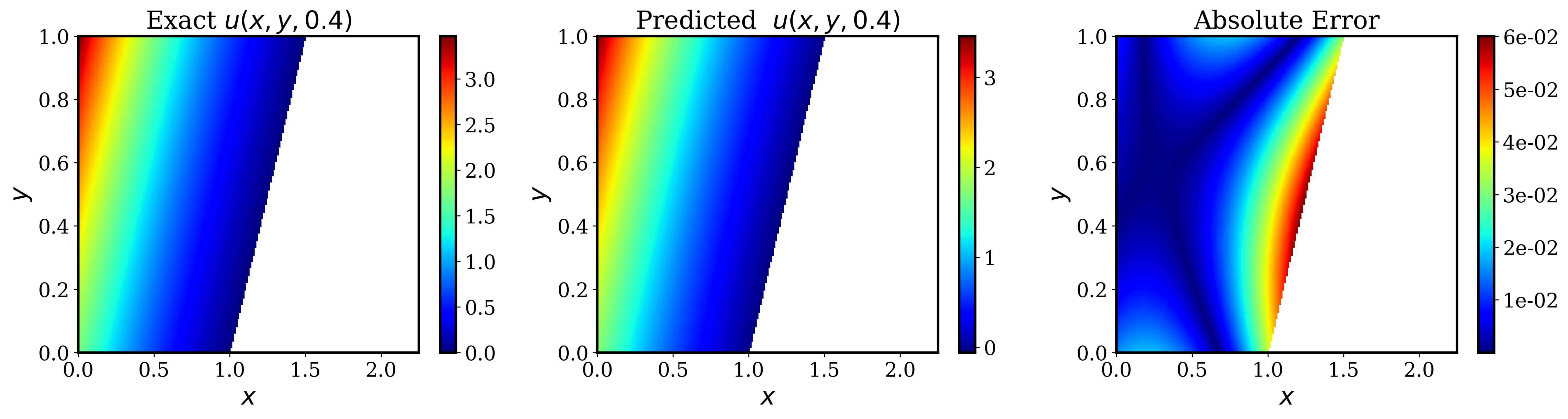}
     \end{subfigure}
     \begin{subfigure}[b]{0.8\textwidth}
         \centering
         \includegraphics[width=\textwidth]{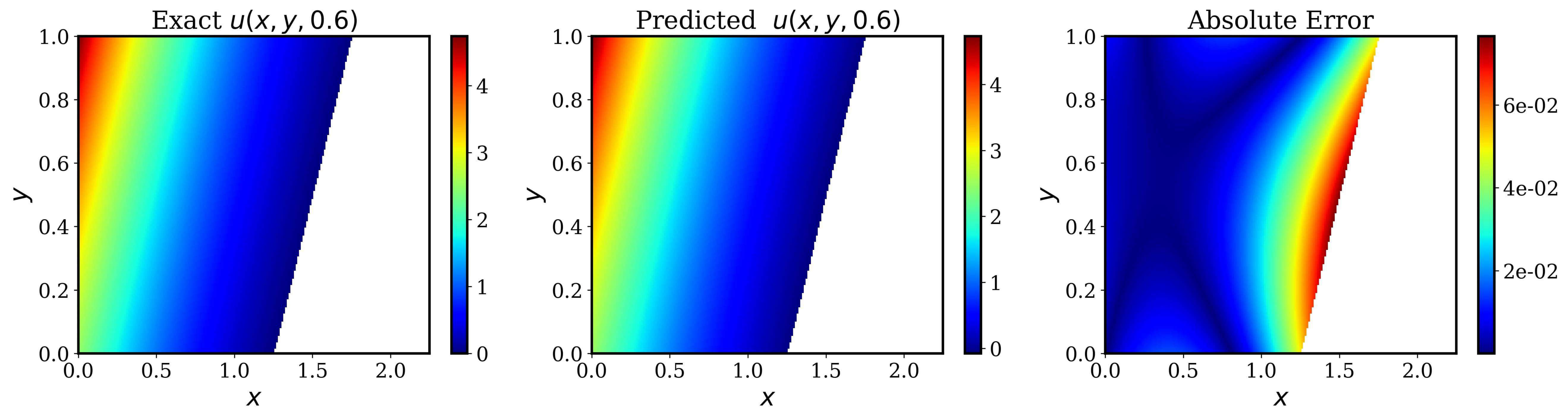}
     \end{subfigure}
      \begin{subfigure}[b]{0.8\textwidth}
         \centering
         \includegraphics[width=\textwidth]{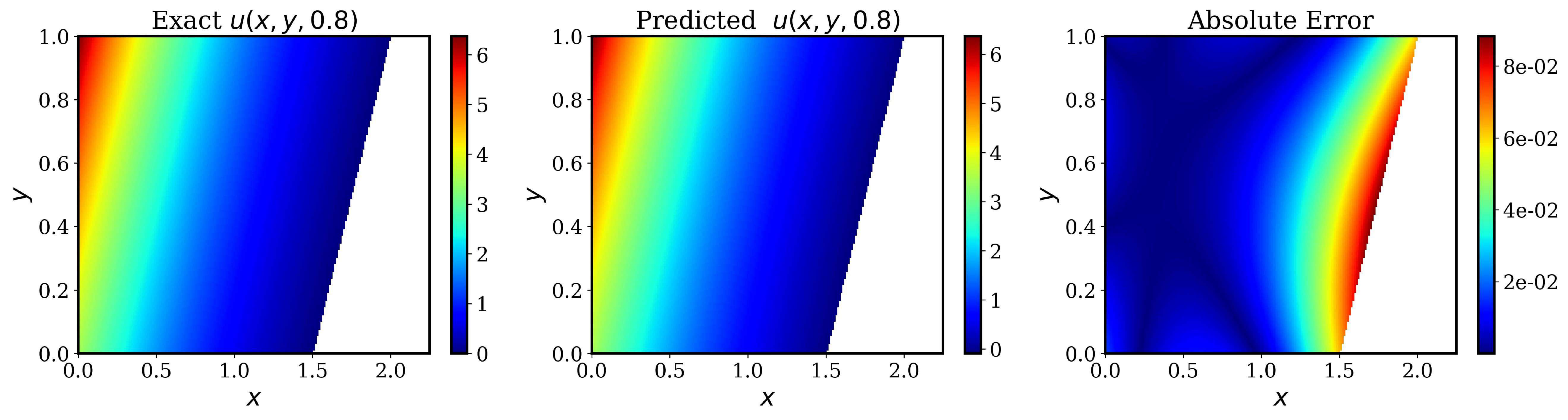}
    
     \end{subfigure}
        \caption{{\em Inverse two-dimensional one-phase Stefan problem Type II:} Comparison of the predicted and exact
solutions corresponding to the four temporal snapshots.}
        \label{fig: Stefan2D_inverse_II_snapshots}
\end{figure}

\begin{figure}
    \centering
    \includegraphics[width=0.8\textwidth]{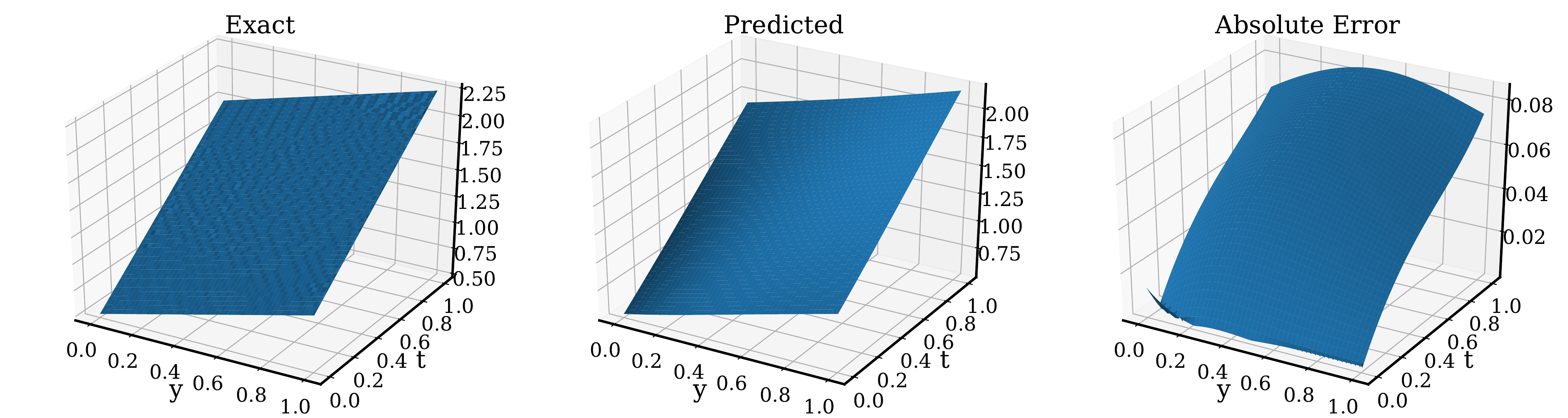}
    \caption{{\em Inverse two-dimensional one-phase Stefan problem Type II:}
    Comparison of the exact and predicted free surface. The relative $L^2$ error is 3.04e-02.}
    \label{fig:Stefan2D_inverse_II_mb}
\end{figure}

We first train our model for the case of noise-free training data. Specifically, the data-set contains $M=50$ randomly sampled sparse measurements of the exact solution. Four temporal snapshots of comparisons between the exact and the predicted solutions are given in figure \ref{fig: Stefan2D_inverse_II_snapshots}, showing a good agreement between the two with an $L^\infty$ error of $O(10^{-2})$. Furthermore, we visualize the identified free surface against the exact one in figure \ref{fig:Stefan2D_inverse_II_mb}, with a relative $L^2$ error of $3.04e-02$. 

To scrutinize the performance of the algorithm further, we have performed a systematic study to quantify its predictive accuracy with respect to the total number of training data, and noise corruption levels $\sigma$.  In particular, we added white noise with magnitude equal to $1\%, 2\%,  5\%, 10\%$ of the $L^\infty$ norm of the solution function $u(x,y,t)$. 
The resulting relative $L^2$ errors for the predicted solution $u(x,y,t)$ and $s(y,t)$ are reported in tables \ref{tab:Stefan_2D1P_inverse_II_pred_u} and table \ref{tab:Stefan_2D1P_inverse_II_pred_mb}, respectively. The general trend shows increased prediction accuracy for both $u(x,y,t)$ and $s(y,t)$ as the total number of training data $M$ is increased, while decreased prediction accuracy as the noise level $\sigma$ is increased. The latter observation comes to no surprise since the data-set becomes less and less accurate due to the higher level of noise. Another crucial observation is that given a sufficient number of training data (e.g., $M=200$), even if the noise corruption is up to $10\%$ (in which case the relative $L^2$ error of predicted temperature solution $u(x,y,t)$ is more than $10\%$), our framework still yields a reasonable identification accuracy for the free boundary $s(y,t)$ with a relative $L^2$ error of $4.48e-02$. Therefore, 
our methodology appears to be quite robust especially for identifying the free boundary with respect to noise levels in the data.

\begin{table}[]
\renewcommand*{\arraystretch}{1.6}
    \centering
    \begin{tabular}{|c|c|c|c|c|c|}
\hline
\diagbox{Data Points M}{Noise $\sigma$}  & $\sigma =0\%$      & $\sigma =1\%$ &  $\sigma =2\%$ &  $\sigma =5\%$&  $\sigma =10\%$ \\ \hline
M = 50 &5.25e-03 &1.85e-02  &6.63e-02 & 1.75e-01   &3.89e-01    \\ \hline
M = 100& 4.65e-03 & 1.93e-02   & 5.09e-02    &6.83e-02   &1.94e-01    \\ \hline
M = 200& 3.72e-03 & 1.31e-02   & 3.48e-02    &8.12e-02 &1.35e-01   \\ \hline
\end{tabular}
    \caption{{\em Inverse two-dimensional one-phase Stefan problem Type II:} Relative $L^2$ errors of predicted solution $u(x,y,t)$ with different number of data points and different noise level $\sigma$.}
    \label{tab:Stefan_2D1P_inverse_II_pred_u}
\end{table}

\begin{table}[]
\renewcommand*{\arraystretch}{1.6}
    \centering
    \begin{tabular}{|c|c|c|c|c|c|}
\hline
\diagbox{Data Points M}{Noise $\sigma$}  & $\sigma =0\%$      & $\sigma =1\%$ &  $\sigma =2\%$ &  $\sigma =5\%$&  $\sigma =10\%$ \\ \hline
M = 50 &3.81e-02 & 3.70e-02    & 4.18e-02  & 1.20e-01   & 1.79e-01    \\ \hline
M = 100&3.64e-02 & 3.98e-02   &  7.54e-02    & 8.13e-02  &  1.51e-01    \\ \hline
M = 200&2.90e-02 & 2.89e-02   &  3.00e-02  & 3.72e-02  &  4.48e-02    \\ \hline
\end{tabular}
    \caption{{\em Inverse two-dimensional one-phase Stefan problem Type II:} Relative $L^2$ errors of predicted free boundary $s(y,t)$ with different number of data points and different noise level $\sigma$.}
    \label{tab:Stefan_2D1P_inverse_II_pred_mb}
\end{table}

\section{Discussion}
\label{sec: summary}
The unique characteristic of free boundary problems is that they pose the challenge of inferring the solution of partial differential equations in domains with unknown boundaries and complex time-dependent interfaces.
Here we have introduced a general deep learning framework for tackling forward and inverse free boundary problems, and tested its effectiveness across a series of numerical case studies involving different formulations of the classical one-phase and two-phase Stefan problems. Moreover, a new type of data-driven inverse Stefan problem has been formulated and addressed under the same unified framework. As demonstrated by the numerical studies presented here, the proposed computational framework is general and flexible in the sense that it requires minimal implementation effort in order to be adapted to different kinds of free boundary problems. 

Although the main focus of this work is to provide a quantitative assessment of how the proposed algorithms perform in a controlled setting, the developments presented here are directly applicable to a wide range of applications in science and engineering involving problems with dynamic free boundaries or phase-transitions. Such example applications include, but are not limited to, fluid-structure interactions \cite{bansch2013ale}, tumor growth modeling \cite{rutter2017mathematical,friedman1999stefan}, thrombus formation \cite{yin2020non}, wound healing \cite{chen2000free}, chemical vapor deposition \cite{friedman2000free}, and electrophotography \cite{chen2000free}. The novel data-driven type of inverse problems put forth in this work is particularly well suited for tackling such applications where a small number of noisy observations may be available.

Despite a series of promising results presented here, there also exist numerous open questions that require further investigation. One of the most challenging questions pertain to handling irregular free boundaries with a complicated geometry that may include sharp cusps, mushy regions, or discontinuities. 
Investigating such cases goes beyond the scope of this work, however we anticipate that the flexibility of deep neural networks in approximating complex functions can provide a new exciting path for addressing these longstanding challenges.

\section*{Acknowledgements}
This work received support from the US Department of Energy under the Advanced Scientific Computing Research program (grant DE-SC0019116) and the Air Force Office of Scientific Research (grant FA9550-20-1-0060).

\bibliographystyle{unsrt}  
\bibliography{references}

\end{document}